# HOW UNIVERSAL ARE ASYMPTOTICS OF DISCONNECTION TIMES IN DISCRETE CYLINDERS?

By Alain-Sol Sznitman

*ETH Zurich*

We investigate the disconnection time of a simple random walk in a discrete cylinder with a large finite connected base. In a recent article of A. Dembo and the author it was found that for large $N$ the disconnection time of $G_N \times \mathbb{Z}$ has rough order $|G_N|^2$, when $G_N = (\mathbb{Z}/N\mathbb{Z})^d$. In agreement with a conjecture by I. Benjamini, we show here that this behavior has broad generality when the bases of the discrete cylinders are large connected graphs of uniformly bounded degree.

**0. Introduction.** We investigate here a simple random walk on an infinite discrete cylinder having its base modeled on a large finite connected graph. We are interested in the time the walk takes to disconnect the cylinder, or in a more picturesque language, in the problem of a "termite in a wooden beam." In a recent work [8], the case when the base is a $d$-dimensional discrete torus of large size $N$, $G_N = (\mathbb{Z}/N\mathbb{Z})^d$, was studied. Answering a question of H. J. Hilhorst, it was shown that for large $N$ the disconnection time typically has rough order $|G_N|^2$. Moreover, it was also conjectured by I. Benjamini that the disconnection time of $G \times \mathbb{Z}$ behaves as $|G|^{2+o(1)}$, for large connected $G$'s of uniformly bounded degree.

We show in this article that the above asymptotic behavior has broad generality and also derive a general asymptotic upper bound on these disconnection times.

We now describe the set-up before discussing the results any further. We consider a finite connected graph with vertex set $G$ and edge set $\mathcal{E}$ made of unordered pairs of $G$. We write $\deg(G)$ for the degree of $G$ (i.e., the maximal number of neighbors of any vertex in $G$). We consider the cylinder based on $G$:

$$(0.1) \qquad E = G \times \mathbb{Z},$$

tacitly endowed with its natural product graph structure. We say that a finite set $S \subseteq E$ disconnects $E$, when for large $M$, $G \times [M, \infty)$ and $G \times (-\infty, -M]$ are contained in two distinct connected components of $E \backslash S$. We denote with $P_x$, $x \in E$, the canonical law on $E^{\mathbb{N}}$ of the simple random walk on $E$ starting at $x$, and with $(X_n)_{n \geq 0}$ the canonical process. We are interested in the disconnection time of $E$:

$$(0.2) \qquad \mathcal{T}_G = \inf\{n \geq 0; X_{[0,n]} \text{ disconnects } E\}.$$

The walk on $E$ is irreducible and recurrent, so that for any $x \in E$, $\mathcal{T}_G$ is $P_x$-a.s. finite. Further, if $\widetilde{\mathcal{C}}_G$ stands for the cover time of $G$ by the projection of $X_.$ on $G$, that is, the first time the projection of $X_.$ has visited all points of $G$, and $\mathcal{C}_G$ stands for the cover time of $G \times \{0\}$ by $X_.$, it is plain that

$$(0.3) \qquad \widetilde{\mathcal{C}}_G \leq \mathcal{T}_G \leq \mathcal{C}_G.$$

There are examples of sequences of finite connected graphs $G_N$ of divergent degree, giving rise to cover times $\widetilde{\mathcal{C}}_{G_N}$ much larger than $|G_N|^2$ (e.g., the "barbells," cf. Aldous and Fill [2], Chapter 5, Example 11). As a result of the left-hand side inequality in (0.3), $\mathcal{T}_{G_N}$ is much larger than $|G_N|^2$ for such sequences. With this in mind, we restrict our attention here to cylinders with bases that are large graphs of uniformly bounded degree. In this context we show in Theorem 1.2 a general upper bound for the disconnection time. Namely, given an integer $d_0$, and $\varepsilon > 0$, one has

$$(0.4) \qquad \lim_{|G| \to \infty, \deg(G) \leq d_0} \sup_{x \in E} P_x[\mathcal{T}_G > |G|^2 (\log|G|)^{4+\varepsilon}] = 0.$$

The above bound exploits the right-hand side inequality in (0.3) and holds as well with $\mathcal{C}_G$ in place of $\mathcal{T}_G$, when the supremum over $E$ in (0.4) is replaced with a supremum over $G \times \{0\}$. We also derive upper bounds on the expectation of $\mathcal{T}_G$ of same order; compare (1.32).

The derivation of a lower bound on $\mathcal{T}_G$ of rough order $|G|^2$ is substantially more delicate. We do not have a lower bound on $\mathcal{T}_G$ of comparable generality to (0.4). The left-hand inequality in (0.3) is now only helpful in a few cases. Indeed, $\widetilde{\mathcal{C}}_G$ is often much smaller than $|G|^2$ (e.g., $\log \widetilde{\mathcal{C}}_{G_N}/\log |G_N|$ is asymptotically close to 1, when $G_N = (\mathbb{Z}/N\mathbb{Z})^d$, with $d \geq 2$, and close to 2 when $d = 1$. For this and much more detailed results, see [1, 5, 7]). In the present work we derive lower bounds on $\mathcal{T}_G$ of "rough order $|G|^2$," when $G$ is large and contains some suitable pocket of possibly vanishing relative volume, inside which we impose additional control. In the pocket we require a quantitative transient or recurrent behavior; see Theorems 4.1 and 5.2. Our methods leave open the case of a too massively recurrent behavior in the pockets; see Theorem 5.2. Otherwise, we also obtain lower bounds on $\mathcal{T}_G$ of rough order $|G|^2$ when the spectral gap $\lambda_G$ [cf. (1.8), (1.9)] is "close" to the extreme possible values compatible with the uniform bound on the degree, that is, for $\lambda_G$ of order $|G|^{-2+o(1)}$ or $|G|^{o(1)}$; see Theorem 4.3.

To give a more explicit flavor of our results, consider, for instance, an infinite connected graph $G_\infty$ of bounded degree. Denoting with $d(\cdot,\cdot)$ the graph distance on $G_\infty$ (i.e., the minimal number of steps of a nearest-neighbor path connecting two points), we assume that for some $\beta \geq 2$, $\alpha \geq (1+\beta/2) \vee (\beta-1)$, and positive constants $\kappa_i, 1 \leq i \leq 4$, one has the sub-Gaussian bounds for the walk $Y$ on $G_\infty$:

$$\text{(i)} \qquad P_g^{G_\infty}[Y_k = g'] \leq \frac{\kappa_1}{k^{\alpha/\beta}} \exp\left\{-\kappa_2 \left(\frac{d(g,g')^\beta}{k}\right)^{1/(\beta-1)}\right\}$$

for $g, g' \in G_\infty$ and $k \geq 1$,

(0.5)

$$\text{(ii)} \quad P_g^{G_\infty}[Y_k \text{ or } Y_{k+1} = g'] \geq \frac{\kappa_3}{k^{\alpha/\beta}} \exp\left\{-\kappa_4 \left(\frac{d(g,g')^\beta}{k}\right)^{1/(\beta-1)}\right\}$$

for $g, g' \in G_\infty$ and $k \geq 1 \vee d(g,g')$.

These bounds are easily seen to imply (cf. Grigoryan and Telcs [10], pages 503–504) that $G_\infty$ is $\alpha$-Ahlfors regular, that is, one has the volume controls

$$(0.6) \qquad \widetilde{\kappa}_1 r^\alpha \leq |B(g,r)| \leq \widetilde{\kappa}_2 r^\alpha \qquad \text{for all } r \geq 1 \text{ and } g \in G_\infty,$$

with $|B(g,r)|$ the cardinality of the open ball in $G_\infty$ with center $g$ and radius $r$, and where the positive numbers $\widetilde{\kappa}_i, i=1,2$, can be chosen as function of $\deg(G_\infty)$ and $\kappa_i, 1 \leq i \leq 4$.

Over the recent years an extensive investigation of such heat-kernel bounds has been made. Equivalent characterizations in terms of volume growth and parabolic Harnack inequality, or mean exit time from balls and Harnack inequality, as well as examples, can be found in Grigoryan and Telcs [10, 11], Barlow [3], Barlow, Coulhon and Kumagai [4] and the references therein. Only values $\alpha + 1 \geq \beta \geq 2$ in (0.5) may and do occur (cf. (2.5) of [10] and [3]) but we are only concerned here with the case $\alpha + 1 \geq (2+\beta/2) \vee \beta$, $\beta \geq 2$ (in particular, this contains the case $\alpha \geq \beta \geq 2$, but excludes certain instances of $\beta > \alpha \geq 1$ that yield so-called very strongly recurrent graphs; cf. [3]). The case $\beta = 2$ was investigated first (cf. Delmotte [6]) and includes usual examples such as $\mathbb{Z}^d$, with $d = \alpha$. The case $\beta > 2$ in (0.5) corresponds to so-called anomalous diffusion, where at time $T$ the walk has traveled at distances of order $T^{1/\beta} \ll \sqrt{T}$; see [3] for examples related to skeletons of fractal sets. We also refer to [11] and [4] for bounds in the context of the more general volume doubling assumption.

As an application of our results, we show in Corollaries 4.5 and 5.3 that when $G_N$ is a sequence of finite connected graphs with cardinality tending to infinity and uniformly bounded degree, such that $G_N$ contains an open ball $A_N$ ("the pocket") isomorphic to some ball in $G_\infty$ satisfying (0.5) (and even less in the transient regime $\alpha > \beta \geq 2$), so that for some sub-polynomially growing sequence $\varphi(n)$ [i.e., $\varphi(n) = o(n^\varepsilon)$, for each $\varepsilon > 0$],

$$(0.7) \qquad |A_N|\varphi(|G_N|) \geq |G_N|,$$

or more generally, such that for some $\eta > 0$ and $\varphi(n)$ as above,

$$(0.8) \qquad \lambda_{G_N}^{1/2}|A_N|\varphi(|G_N|) \geq \min(|G_N|^\eta, \lambda_{G_N}^{1/2}|G_N|) \qquad \text{for large } N,$$

then for any $\delta, \varepsilon > 0$, writing $E_N = G_N \times \mathbb{Z}$, one has

$$(0.9) \qquad \lim_{N \to \infty} \inf_{x \in E_N} P_x[|G_N|^{2(1-\delta)} \leq \mathcal{T}_{G_N} \leq |G_N|^2 (\log|G_N|)^{4+\varepsilon}] = 1.$$

This result alone covers many examples and vastly generalizes Theorems 1.1 and 2.1 of [8]. Our methods, however, leave open the case $\beta \geq 2$, and $1 \leq \alpha < 1 + \beta/2$, corresponding to some instances of so-called very strong recurrence of $G_\infty$; see [3]. We otherwise have applications beyond the above set-up. For instance, we show in Corollary 4.5 that (0.9) holds true when $G_N$ is the rooted $r$-tree of depth $N$, or also when (cf. Remark 4.4),

$$(0.10) \qquad \lambda_{G_N} = |G_N|^{o(1)} \quad \text{or} \quad \lambda_{G_N} = |G_N|^{-2+o(1)}.$$

We now give some indications on the techniques we employ in this work. As already mentioned, lower bounds on the disconnection time cause the main difficulty. The strategy in this work differs in several respects from the line followed in [8], when $G$ is the $d$-dimensional torus of size $N$. In [8] a crucial role was played by the geometric Lemmas 2.4 and 2.5, which show that when $S$ disconnects $E$, one can find on a whole range of scales cubes in $E$ where $S$ has a trace with cardinality, which is at least that of a fraction of a face of the cube. The length scale is then adjusted so that typically up to time $|G|^{2(1-\delta)}$, for any cube of corresponding side-length, few excursions of the walk enter the cube, and the walk can hardly leave a trace comparable in cardinality to the face of the cube. Implicit to this approach are certain isoperimetric controls that need not hold true in our context. To give a feel for the issue, observe that in a rooted binary tree of finite depth, unlike what happens for discrete tori of dimension $d \geq 2$, one can find subsets of roughly half volume with boundary consisting of a single point (the root). Thus, insisting on isoperimetric controls of the type used in [8] rules out many interesting examples.

We follow here a different route. We construct with high probability connections between top and bottom of the cylinder that avoid the trajectory $X$ up to time $|G|^{2(1-\delta)}$. We use a localization technique that enables to focus on what happens in a sub-cylinder $A \times \mathbb{Z}$ of $E$, with $A$ the "pocket," a possibly very small subset of $G$. We analyze excursions of the walk entering a suitably small box $C$, with base sitting well inside the pocket $A$ [cf. (3.11)], which then move at vertical distances of order $2h'$ from $C$; see (3.6). We show in Proposition 3.2 that typically only finitely many such excursions occur up to time $|G|^{2(1-\delta)}$. The height $h'$ is, on the one hand, chosen big enough so that starting from a point with $G$-projection inside $A$, at vertical distance of order $h'$ from $C$, the walk has a small enough probability of entering $C$ before moving at vertical distance $2h'$ from $C$. On the other hand, $h'$ is chosen sufficiently small so that what happens outside $A \times \mathbb{Z}$ has little influence on what happens inside $C$. As a by-product, the finitely many excursions that typically enter $C$ are also of truly shorter duration than the naive excursions employed in Section 1, for which $h'$ is replaced with a height $h$ slightly bigger than $\lambda_G^{-1/2}$; compare (1.17) and (1.10). This makes it easier to control the damage they may cause inside $C$.

Rarefaction of excursions to $C$ is the first step in constructing many top-to-bottom connections in a sub-box $D$ of $C$, which avoid the walk up to time $|G|^{2(1-\delta)}$. The second step consists in containing the damage the finitely many excursions reaching $C$ may create. We rely here on ensuring sufficiently many horizontal and vertical connections across certain boxes, and a renormalization procedure, which is used when the walk has recurrent behavior in the pocket. In this fashion we construct with high probability very connective boxes $D$ that can be piled up to produce top-to-bottom connections in the cylinder $E$; see Proposition 2.6. As already hinted at, handling recurrent pockets in $G$ is more delicate than dealing with transient pockets, and leads us to require additional control; see (5.2) and (5.3).

We now describe the organization of the article.

In Section 1 we introduce additional notation and mainly derive the general upper bound (0.4) on $\mathcal{T}_G$ in Theorem 1.2. The essential point is to bound the cover time of $G \times \{0\}$ from above.

In Section 2 we develop auxiliary results that are preparatory for the lower bound on $\mathcal{T}_G$. These results pertain to the localization technique (cf. Proposition 2.3) to the construction of connective blocks [cf. (2.43) and Proposition 2.6] and to the treatment of graphs with low lying spectral gap, see Proposition 2.1.

In Section 3 we develop the localization technique and show in Proposition 3.2 that few excursions of the walk meet the box $C$ by time $|G|^{2(1-\delta)}$.

In Section 4 we derive a lower bound on $\mathcal{T}_G$ in the case of a transient pocket (cf. Theorem 4.1), or when the spectral gap is close to its extreme values; see Theorem 4.3. Applications are given in Corollaries 4.5, 4.6 and Remark 4.7.

In Section 5 we obtain a lower bound on $\mathcal{T}_G$ that applies to cases of recurrent behavior in the pocket; see Theorem 5.2. Applications are then given in Corollary 5.3 and Remark 5.5.

**1. The upper bound.** The main object of this section is to prove a general asymptotic upper bound on the disconnection time of discrete cylinders based on large finite connected graphs of uniformly bounded degree. The principal result appears in Theorem 1.2, where, in particular, (0.4) is derived. The proof exploits the right-hand side inequality of (0.3) and mainly focuses on bounding the cover time $\mathcal{C}_G$ of $G \times \{0\}$ from above. We first introduce additional notation, and recall some classical facts.

For $u$ a nonnegative real number, we let $[u]$ stand for the integer part of $u$. For $v, w$ real numbers, we write $v \wedge w$ and $v \vee w$ for the minimum and the maximum of $v$ and $w$. Given a finite set $A$, we denote with $|A|$ its cardinality. When $\Gamma$ is a graph and $x, x'$, are distinct vertices of $\Gamma$, we write $x \sim x'$, if $x$ and $x'$ are neighbors, that is, $\{x, x'\}$ is an edge of $\Gamma$; we denote with $\deg(x)$ [or $\deg_\Gamma(x)$ if there is a risk of confusion] the degree of $x$, that is, the number of neighbors of $x$, and $\deg(\Gamma) = \sup\{\deg(x); x \text{ vertex of } \Gamma\}$, the degree of $\Gamma$. With an abuse of notation we usually make no distinction between a graph and its set of vertices. We denote with $d(\cdot, \cdot)$ [or sometimes with $d_\Gamma(\cdot, \cdot)$], the distance function on $\Gamma$, that is, the minimal number of steps for a nearest neighbor path on $\Gamma$ joining two given points of $\Gamma$. The graphs we consider in the sequel are all connected so that $d(\cdot, \cdot)$ is automatically finite. We denote with $B(x, r)$ [or $B_\Gamma(x, r)$ when there is a risk of confusion] the open ball with center $x \in \Gamma$ and radius $r > 0$. When $U$ is a subset of $\Gamma$, we denote with $\partial U$ its boundary:

$$(1.1) \qquad \partial U = \{x \in U^c; \exists x' \in U \text{ with } x \sim x'\}.$$

Throughout the article the finite connected graphs $G$ (with edge set $\mathcal{E}$) that show up as the base of the cylinder $E = G \times \mathbb{Z}$ have degree uniformly bounded by some integer $d_0 \geq 2$,

$$(1.2) \qquad \deg(G) \leq d_0 \text{ and we tacitly assume } |G| \geq 2.$$

Since $G$ is connected, it follows that

$$(1.3) \qquad |G| \leq 2|\mathcal{E}| \leq d_0 |G|.$$

We write $\pi_G$ and $\pi_\mathbb{Z}$ for the respective canonical projections of $E$ on $G$ and $\mathbb{Z}$.

We denote with $X., Y., Z.$ the respective canonical walks in discrete time on $E, G, \mathbb{Z}$, which at each step jump with equal probability to one of the neighbors of their current location. We write $P_x, P_g^G, P_u^\mathbb{Z}$ for the respective canonical laws starting at $x \in E$, $g \in G$, $u \in \mathbb{Z}$. The canonical shifts and filtrations are denoted with $(\theta_n)_{n \geq 0}$ and $(\mathcal{F}_n)_{n \geq 0}$, with a possible superscript $E$, $G$ or $\mathbb{Z}$, when confusion may arise. For a subset $U$ of $E$, $G$ or $\mathbb{Z}$, we denote with $H_U$ and $T_U$ the entrance time in $U$ and exit time from $U$ of the respective walk, so, for instance, when $U \subseteq E$,

$$(1.4) \qquad H_U = \inf\{n \geq 0, X_n \in U\}, \qquad T_U = \inf\{n \geq 0, X_n \notin U\},$$

with $X.$ replaced by $Y.$ or $Z.$, when $E$ is replaced by $G$ or $\mathbb{Z}$. Again, when confusion may arise, we add a superscript $G$ or $\mathbb{Z}$ to clarify the notation. When $U$ is a singleton $\{z\}$, we write $H_z$ in place of $H_{\{z\}}$.

It is convenient to consider the canonical continuous time random walks $\overline{X}., \overline{Y}., \overline{Z}.$, which respectively jump with rates $\deg(g)+2$, $\deg(g)$ and $2$, when respectively located at $x = (g,u)$, $g$ and $u$. With an abuse of notation, we still denote with $P_x$, $P_g^G$, $P_u^{\mathbb{Z}}$ the corresponding canonical laws. Otherwise, we use notation such as $(\overline{\theta}_t)_{t \geq 0}$, $(\overline{\mathcal{F}}_t)_{t \geq 0}$ or $\overline{H}_U$ to refer to the natural continuous time objects. Clearly, the respective discrete skeletons of the continuous time walks $\overline{X}., \overline{Y}., \overline{Z}.$ are distributed as the respective discrete time walks $X., Y., Z.$. Further, the continuous time walks satisfy the following useful fact, that we recurrently use in the sequel:

$$\text{(1.5)} \quad \begin{array}{l} \text{for } x = (g, u) \in E, \quad \text{under } P_g^G \otimes P_u^{\mathbb{Z}}, \\ (\overline{Y}., \overline{Z}.) \text{ has the canonical law } P_x \text{ governing } \overline{X}.. \end{array}$$

The stationary distributions of the discrete and continuous time walks on $G$ are the reversible measures (for the respective walks) defined by

$$\text{(1.6)} \quad \mu(g) = \frac{\deg(g)}{2|\mathcal{E}|}, \qquad \overline{\mu}(g) = \frac{1}{|G|} \qquad \text{for } g \in G.$$

The generator and the Dirichlet form attached to the continuous time walk on $G$ are respectively

$$\text{(1.7)} \quad \begin{array}{l} L_G f(g) = \sum_{g' \sim g} (f(g') - f(g)), \qquad g \in G, \\[2mm] \mathcal{D}_G(f, f) = (-L_G f, f)_{L^2(\overline{\mu})} = \dfrac{1}{2|G|} \sum_{\substack{g, g' \in G \\ g \sim g'}} (f(g') - f(g))^2, \end{array}$$

with $f$ an arbitrary function on $G$, and $(\cdot, \cdot)_{L^2(\overline{\mu})}$ the $L^2$-scalar product on $G$. In what follows an important role is played by the spectral gap of the continuous time walk on $G$:

$$\text{(1.8)} \quad \lambda_G = \inf_{f \text{ nonconstant}} \frac{\mathcal{D}_G(f, f)}{\mathrm{var}_{\overline{\mu}}(f)} \quad \text{with } \mathrm{var}_{\overline{\mu}}(f) \text{ the variance of } f \text{ under } \overline{\mu}.$$

It follows from Cheeger's inequality (cf. Aldous and Fill [2], Chapter 4, Section 5.2, page 34, or Lubotzky [13], Propositions 4.2.4 and 4.2.5), for the lower bound and the choice in (1.8) of a function $f$ vanishing everywhere except at a single point of $G$, that

$$\text{(1.9)} \quad 2d_0 \geq \lambda_G \geq \frac{2}{d_0 |G|^2}.$$

We introduce the time

(1.10) $$t_G = \lambda_G^{-1} \log(2|G|),$$

which will play an important role in the sequel, due to the following (classical) result:

LEMMA 1.1.

(1.11)
$$\text{For } t \geq t_G, g, g' \in G,$$
$$|(P_g^G[\overline{Y}_t = g']/\overline{\mu}(g')) - 1| \leq \tfrac{1}{2} \exp\{-(t - t_G)\lambda_G\}.$$

PROOF. The argument is classical; see Saloff-Coste [14], page 328. Writing

(1.12) $$p_t(g, g') = P_g^G[\overline{Y}_t = g']\overline{\mu}(g')^{-1}, \qquad g, g' \in G, t \geq 0,$$

for the transition density of the continuous walk on $G$, it follows from the spectral theorem that

$$\sum_{g' \in G} (p_t(g, g') - 1)^2 \overline{\mu}(g') \leq e^{-2\lambda_G t} \sum_{g' \in G} (1_{\{g'=g\}} \overline{\mu}(g)^{-1} - 1)^2 \overline{\mu}(g')$$
$$= e^{-2\lambda_G t}(\overline{\mu}(g)^{-1} - 1).$$

The claim (1.11) then follows from the fact that $\exp\{-2\lambda_G t_G\} = (4|G|^2)^{-1}$, and (1.6). $\square$

We now introduce certain stopping times that will be used throughout the article. Given an integer $h \geq 1$ and $u \in \mathbb{Z}$, we consider the boxes in $E$:

(1.13)
$$B_h(u) = G \times I(u) \subseteq \widetilde{B}_h(u) = G \times \widetilde{I}(u)$$
$$\text{with } I(u) = u + [-h, h] \quad \text{and} \quad \widetilde{I}(u) = u + [-2h + 1, 2h - 1].$$

We write $B_h$, $\widetilde{B}_h$ in place of $B_h(0)$, $\widetilde{B}_h(0)$, and when the value of $h$ is clearly specified, we simply drop the subscript $h$ from the notation. The successive returns to $B_h(u)$ and departures of $\widetilde{B}_h(u)$ are then defined by

(1.14)
$$R_1^{h,u} = H_{B_h(u)},$$
$$D_1^{h,u} = T_{\widetilde{B}_h(u)} \circ \theta_{R_1^{h,u}} + R_1^{h,u}, \quad \text{and for } k \geq 1,$$
$$R_{k+1}^{h,u} = R_1^{h,u} \circ \theta_{D_k^{h,u}} + D_k^{h,u}, \qquad D_{k+1}^{h,u} = D_1^{h,u} \circ \theta_{D_k^{h,u}} + D_k^{h,u},$$

so that

$$0 \leq R_1^{h,u} \leq D_1^{h,u} \leq \cdots \leq R_k^{h,u} \leq D_k^{h,u} \leq \cdots \leq \infty,$$

and for any $x \in E$, $P_x$-a.s., these inequalities are strict except maybe the first one. With a similar convention as above, we drop the superscript $h$ when the value of $h$ is clearly specified and the superscript $u$ when $u = 0$.

Let us explain our convention concerning constants for the remainder of this section and Section 2 as well. We will denote with $c$ a positive constant solely depending on $d_0$ [cf. (1.2)], with value changing from place to place. Additional dependence will appear in the notation, for instance, $c(\varepsilon)$ refers to a positive constant depending on $d_0$ and $\varepsilon$. Numbered constants like $c_0, c_1, \ldots$ will refer to the value of the constant in the first display where they are determined. Finally, we will use the expression for large $G$, in place of for $|G| \geq c$, with $G$ a finite connected graph satisfying (1.2). The main result of this section is the following:

THEOREM 1.2.

$$(1.15) \quad \lim_{|G| \to \infty, \deg(G) \leq d_0} \inf_{x \in G \times \{0\}} P_x[\mathcal{T}_G \leq \mathcal{C}_G \leq |G|^2 (\log |G|)^{4+\varepsilon}] = 1$$

*for any $\varepsilon > 0$.*

REMARK 1.3. It is plain that, for any $g \in G, u \in \mathbb{Z}$, the disconnection time $\mathcal{T}_G$ has the same distribution under $P_{(g,u)}$ and $P_{(g,0)}$, so that (1.15) readily implies

$$(1.16) \quad \lim_{|G| \to \infty, \deg(G) \leq d_0} \inf_{x \in E} P_x[\mathcal{T}_G \leq |G|^2 (\log |G|)^{4+\varepsilon}] = 1 \qquad \text{for any } \varepsilon > 0.$$

PROOF OF THEOREM 1.2. Throughout the remainder of this section the value of $h$ [cf. (1.13), (1.14)] is set equal to

$$(1.17) \qquad h = [\sqrt{t_G}] + 2.$$

For any $z = (g', 0) \in G \times \{0\}$ and $x = (g, u) \in B$ [cf. (1.13) below, for the notation] the strong Markov property for the continuous time walk at time $\overline{H}_z$ implies that

$$(1.18) \qquad P_x[H_z < T_{\widetilde{B}}] = P_x[\overline{H}_z < \overline{T}_{\widetilde{B}}] = \frac{a_1}{a_2},$$

where

$$(1.19) \quad \begin{aligned} a_1 &= E_x\left[\int_0^\infty 1\{\overline{X}_t = z, t < \overline{T}_{\widetilde{B}}\}\, dt\right], \\ a_2 &= E_z\left[\int_0^\infty 1\{\overline{X}_t = z, t < \overline{T}_{\widetilde{B}}\}\, dt\right]. \end{aligned}$$

We now bound $a_1$ from below and $a_2$ from above, and thus obtain a lower bound on the left-hand side member of (1.18). With the help of (1.5) and the notation (1.13), we find that

$$a_1 = \int_0^\infty P_g^G[\overline{Y}_t = g'] P_u^{\mathbb{Z}}[\overline{Z}_t = 0, t < \overline{T}_{\widetilde{I}}] \, dt$$

$$(1.20) \quad \overset{(1.11)}{\geq} \int_{t_G}^\infty \frac{1}{2|G|} P_u^{\mathbb{Z}}[\overline{Z}_t = 0, t < \overline{T}_{\widetilde{I}}] \, dt$$

$$\geq \frac{1}{2|G|} E_u^{\mathbb{Z}}[t_G < \overline{T}_{\widetilde{I}}, P_{\overline{Z}_{t_G}}[\overline{H}_0 < \overline{T}_{\widetilde{I}}]] E_0^{\mathbb{Z}}\left[\int_0^{\overline{T}_{\widetilde{I}}} 1\{\overline{Z}_t = 0\} \, dt\right],$$

using the strong and the simple Markov property in the last step. It follows from the invariance principle and (1.17) that the first expectation in the last line of (1.20) is bounded below by a positive constant. Using standard calculations on the continuous and discrete simple random walk on $\mathbb{Z}$, we also find that

$$E_0^{\mathbb{Z}}\left[\int_0^{\overline{T}_{\widetilde{I}}} 1\{\overline{Z}_t = 0\} \, dt\right] = \tfrac{1}{2} E_0^{\mathbb{Z}}\left[\sum_{k \geq 0} 1\{Z_k = 0, k < T_{\widetilde{I}}\}\right] \geq ch.$$

Collecting the lower bounds we have derived, we find that

$$(1.21) \qquad a_1 \geq \frac{ch}{|G|} \quad \text{[with } h \text{ defined by (1.17)]}.$$

We will now obtain an upper bound on $a_2$ in (1.19). We first note that

$$(1.22) \qquad |P_g^G[\overline{Y}_t = g'] - \overline{\mu}(g')| \leq \frac{c}{\sqrt{t}} \qquad \text{for } t > 0, g, g' \in G$$

(see the convention concerning constants stated above Theorem 1.2). Indeed, (1.22) follows from Theorem 2.3.1, page 345 of Saloff-Coste [14], and the Nash-type inequality

$$(1.23) \qquad \operatorname{var}_{\overline{\mu}}(f)^3 \leq c|G|^2 \mathcal{D}_G(f,f) \|f\|_{L^1(\overline{\mu})}^4$$

for $f$ an arbitrary function on $G$.

The above inequality (1.23) is proven in the same fashion as described in Example 2.3.1, pages 348–350, of [14]. For a related inequality to (1.22), we also refer to Proposition 18 in Chapter 6, Section 4.2 of Aldous and Fill [2]. Therefore, in view of (1.19) [recall $z = (g', 0) \in G \times \{0\}$], we find

$$a_2 = \int_0^\infty P_z[\overline{X}_t = z, t < \overline{T}_{\widetilde{B}}] \, dt$$

$$= \sum_{k \geq 0} \int_{kt_G}^{(k+1)t_G} P_z[\overline{X}_t = z, t < \overline{T}_{\widetilde{B}}] \, dt$$

(1.24)

$$\leq \sum_{k\geq 0} P_z[kt_G < \overline{T}_{\widetilde{B}}] \int_0^{t_G} P_z[\overline{X}_t = z]dt$$

$$\stackrel{(1.5)}{=} \sum_{k\geq 0} P_0^{\mathbb{Z}}[kt_G < \overline{T}_{\widetilde{I}}] \int_0^{t_G} P_z[\overline{X}_t = z]dt,$$

where in the second line we have used the simple Markov property at time $kt_G$ followed by the strong Markov property at time $\overline{H}_z$. From the invariance principle and the Markov property at times $\ell t_G$, $0 \leq \ell < k$, we infer that

(1.25) $\qquad P_0^{\mathbb{Z}}[kt_G < \overline{T}_{\widetilde{I}}] \leq e^{-ck} \qquad$ for any $k \geq 0$.

Coming back to (1.24), we thus find that

$$a_2 \quad \leq \quad c\int_0^{t_G} P_z[\overline{X}_t = z]\,dt \stackrel{(1.5)}{=} c\int_0^{t_G} P_{g'}^G[\overline{Y}_t = g']P_0^{\mathbb{Z}}[\overline{Z}_t = 0]\,dt$$

(1.26) $\qquad \stackrel{(1.22)}{\leq} \quad c\int_0^{t_G} \left[\left(\frac{c}{\sqrt{t}} + \frac{1}{|G|}\right) \wedge 1\right]\frac{1}{\sqrt{t}}\,dt \leq c\frac{\sqrt{t_G}}{|G|} + c\log t_G$

$$\stackrel{(1.9),(1.10)}{=} c\log t_G \leq c\log |G|.$$

Coming back to (1.18), it thus follows from (1.21) and (1.26) that

(1.27) $\quad P_x[H_z < T_{\widetilde{B}}] \geq \dfrac{ch}{|G|\log |G|} \qquad$ for any $z \in G \times \{0\}$ and $x \in B$.

If we now apply the strong Markov property at times $R_m$, $m \geq 1$ [cf. (1.14) and below (1.14) for the notation] we thus find that, for $k \geq 1$, $x \in B$, $z \in G \times \{0\}$,

(1.28) $\quad P_x[H_z > R_k] \leq \left(1 - \dfrac{ch}{|G|\log |G|}\right)^{k-1} \leq \exp\left\{-c_1 \dfrac{h(k-1)}{|G|\log |G|}\right\}.$

We then set $c_2 = 2c_1^{-1}$ and define

(1.29) $\qquad k_* = \left[c_2 \dfrac{|G|}{h}(\log |G|)^2\right] + 2.$

Note that in view of (1.9), (1.10) and (1.17), $\lim_{|G|\to\infty, \deg(G)\leq d_0} k_* = \infty$. We now see that for $x \in B$,

$$P_x[\mathcal{C}_G > R_{k_*}] \leq \sum_{z\in G\times\{0\}} P_x[H_z > R_{k_*}] \stackrel{(1.28)}{\leq} |G|\exp\left\{-c_1\frac{h(k_*-1)}{|G|\log |G|}\right\}$$
$$\leq |G|\exp\{-c_1 c_2(\log |G|)\} = \frac{1}{|G|}.$$

We have thus obtained that

(1.30) $\qquad \lim_{|G|\to\infty, \deg(G)\leq d_0} \sup_{x\in B} P_x[\mathcal{C}_G > R_{k_*}] = 0.$

With similar bounds as (1.21) and (1.22) in [8] [see also (1.18) of the same reference] we can bound $R_{k_*}$, and find that, for any $\varepsilon > 0$,

$$(1.31) \qquad \lim_{|G|\to\infty, \deg G \leq d_0} \sup_{x \in B} P_x[R_{k_*} > (k_* h)^2 (\log|G|)^\varepsilon] = 0.$$

With (1.29), (1.30) and (0.3), this is more than enough to prove (1.15). Incidentally, let us mention that there is some flexibility with the choice of $h$ in (1.17), and the above proof works with minor changes in (1.24)–(1.26), if for large $G$ we choose $h$ as a positive integer lying between $\sqrt{t_G}$ and $|G|\log|G|$. □

REMARK 1.4. Theorem 1.2 also leads to our upper bound on $\sup_{x \in E} E_x[\mathcal{T}_G]$. Indeed, it follows from Theorem 1.2 that, for any $\varepsilon > 0$, when $|G| \geq c(\varepsilon)$,

$$\inf_{x \in E} P_x[\mathcal{T}_G < |G|^2 (\log|G|)^{4+\varepsilon/2}] \geq \tfrac{1}{2},$$

so that with simple Markov property and $\mathcal{W}_G = \mathcal{T}_G/(|G|^2 (\log|G|)^{4+\varepsilon/2})$,

$$\sup_{x \in E} P_x[\mathcal{W}_G \geq k] \leq (\tfrac{1}{2})^k \qquad \text{for } k \geq 0, \text{ whence } \sup_{x \in E} E_x[\mathcal{W}_G] \leq 2.$$

We thus find that, for any $\varepsilon > 0$,

$$(1.32) \qquad \lim_{|G|\to\infty, \deg G \leq d_0} \sup_{x \in E} \frac{E_x[\mathcal{T}_G]}{|G|^2 (\log|G|)^{4+\varepsilon}} = 0.$$

Incidentally, note that in contrast to (1.32), due to the nonintegrability of the hitting time (i.e., first entrance time after time 1) of 0, for the simple random walk on $\mathbb{Z}$, $E_x[\mathcal{C}_G] = \infty$ for all $x \in E$.

**2. Some auxiliary results.** In this section we discuss four auxiliary results that will be helpful in the derivation of lower bounds on the disconnection times of discrete cylinders in the next two sections. The first result (see Proposition 2.1) shows in a quantitative way that the disconnection of $E$ typically cannot take place up to times almost of order $\lambda_G^{-1}$. The next result (cf. Proposition 2.3) is part of the localization technique that enables to focus on what happens in the sub-cylinder $A \times \mathbb{Z}$ of $E$, when $A \subseteq G$ is suitably chosen. The third result (cf. Lemma 2.5) provides upper bounds on the probability that the walk hits a point before exiting $\widetilde{B}$ [cf. (1.14) and (2.12)] and yields exponential controls on the $G$- and $\mathbb{Z}$-projections of the trace in a sub-cylinder of $E$ of the trajectory of the walk up to the time it exits $\widetilde{B}$. These controls will especially be helpful in Section 3 to handle the case of "high values" of $\lambda_G$. The fourth result (cf. Proposition 2.6) describes the basic strategy we employ, when proving that disconnection of the cylinder does not take place up to a certain time. In some sense it replaces and by-passes the arguments based on isoperimetric controls that were used in [8] (cf. Lemmas 2.4 and 2.5) in the case of $G = (\mathbb{Z}/N\mathbb{Z})^d$, with $d \geq 2$. Throughout this section we keep the same convention concerning constant and the use of the expression "for large $G$" as explained above Theorem 1.2.

We first introduce some additional notation. The kernel of the simple random walk on $G$ is

$$R_G f(g) = \deg(g)^{-1} \sum_{g' \sim g} f(g') \tag{2.1}$$

for $g \in G$ and $f$ an arbitrary function on $G$.

We consider $\varphi$ a normalized eigenfunction of $L_G$ [cf. (1.7)], attached to $-\lambda_G$ [cf. (1.8)]

$$-L_G \varphi = \lambda_G \varphi \quad \text{with } \sum_{g \in G} \varphi^2(g) \overline{\mu}(g) = 1, \text{ which then automatically} \tag{2.2}$$

satisfies the orthogonality condition $\sum_{g \in G} \varphi(g) \overline{\mu}(g) = 0$.

We also denote with $W$ the subset of $G$:

$$W = \{g \in G; \varphi(g) > 0\} \subseteq G. \tag{2.3}$$

The first result of this section is the following:

PROPOSITION 2.1.

$$\text{For } n \geq 0, g \in W, \quad P_g^G[T_W > n] \geq \frac{\varphi(g)}{\max \varphi}(1 - \lambda_G)_+^n, \tag{2.4}$$

$$\lim_{|G| \to \infty, \deg G \leq d_0} \sup_{x \in E} P_x[\mathcal{T}_G \leq \lambda_G^{-1} \varepsilon_{|G|}] = 0 \tag{2.5}$$

for any positive sequence $\varepsilon_n$ with $\lim_n \varepsilon_n = 0$.

PROOF.  We begin with the proof of (2.4). We need only consider the case $\lambda_G < 1$. From (1.7), (2.1) and (2.2), we find that

$$(R_G \varphi)(g) = \left(1 - \frac{\lambda_G}{\deg(g)}\right) \varphi(g) \quad \text{for } g \in G. \tag{2.6}$$

As a result, we see that

$$\varphi \text{ and } R_G \varphi \text{ are positive} \quad \text{on } W, \text{ and}$$
$$0 < \frac{\varphi}{R_G \varphi} \leq (1 - \lambda_G)^{-1} \quad \text{on } W. \tag{2.7}$$

Writing $\varphi_+$ for $\max(\varphi, 0)$ and applying the stopping theorem to the $(\mathcal{F}_n^Y)$-martingale

$$\psi_\varepsilon(Y_n) \prod_{k=0}^{n-1} \left(\frac{\psi_\varepsilon}{R_G \psi_\varepsilon}\right)(Y_k), \qquad n \geq 0, \text{ where } \psi_\varepsilon = \varphi_+ + \varepsilon, \text{ with } \varepsilon > 0,$$

we see, using dominated convergence, (2.7), and letting $\varepsilon$ tend to 0, that for $n \geq 0$, $g \in W$,

$$
\begin{aligned}
\varphi(g) &= E_g^G\left[\varphi_+(Y_{n \wedge T_W}) \prod_{k=0}^{n \wedge T_W - 1} \frac{\varphi}{R_G \varphi_+}(Y_k)\right] \\
&= E^G\left[\varphi(Y_n) \prod_{k=0}^{n-1} \frac{\varphi}{R_G \varphi_+}(Y_k), n < T_W\right] \\
&\quad + E_g^G\left[\varphi_+(Y_{T_W}) \prod_{k=0}^{T_W - 1} \frac{\varphi}{R_G \varphi_+}(Y_k), n \geq T_W\right] \\
&\leq E^G\left[\varphi(Y_n) \prod_{k=0}^{n-1} \frac{\varphi}{R_G \varphi}(Y_k), n < T_W\right] \\
&\stackrel{(2.7)}{\leq} \max \varphi (1 - \lambda_G)^{-n} P_g^G[n < T_W],
\end{aligned}
$$
(2.8)

since the first term in the second line vanishes and $R_G \varphi_+ \geq R_G \varphi > 0$, on $W$. The claim (2.4) follows.

We then turn to the proof of (2.5). Without loss of generality, we assume that $\lambda_G < \frac{1}{2}$ [indeed, in the case of graphs with $\lambda \geq \frac{1}{2}$, (2.5) becomes obvious]. With $\varphi$ as above, we pick $g_+ \in G$ such that $\varphi(g_+) = \max \varphi$. It then follows from (2.4) that, for $n \geq 0$,

$$P_{g_+}^G[T_W \leq n] \leq 1 - (1 - \lambda_G)^n \leq n \log\left(\frac{1}{1 - \lambda_G}\right). \tag{2.9}$$

A similar inequality holds for the set $V = \{\varphi < 0\}$ and $g_- \in V$ such that $\varphi(g_-) = \min \varphi$, in place of $W$ and $g_+$, respectively. We then introduce the $(\mathcal{F}_n^Y)$-stopping times:

$$
\begin{aligned}
\tau &= \inf\{n \geq 1; \varphi(Y_n)\varphi(Y_{n-1}) \leq 0\} \quad \text{and} \\
\rho &= \tau \circ \theta_{H_{\{g_+, g_-\}}} + H_{\{g_+, g_-\}},
\end{aligned}
$$
(2.10)

in other words, $\rho$ is the first time $\varphi(Y_n)$ changes sign after reaching either $g_-$ or $g_+$. Given any sequence $\varepsilon_n$ as in (2.5), one has for any $x = (g, u) \in E$,

$$P_x[\mathcal{T}_G \leq \lambda_G^{-1} \varepsilon_{|G|}] \leq P_g^G[\rho \leq \lambda_G^{-1} \varepsilon_{|G|}]$$

$$(2.11) \qquad \leq E_g^G[P_{Y_{H_{\{g_+,g_-\}}}}^G[\tau \leq \lambda_G^{-1}\varepsilon_{|G|}]]$$

$$\stackrel{(2.9)}{\leq} \lambda_G^{-1}\varepsilon_{|G|}\log\left(\frac{1}{1-\lambda_G}\right).$$

Observing that the function $s \in (0, \frac{1}{2}] \to s^{-1}\log(\frac{1}{1-s})$ is bounded, (2.5) follows. $\square$

REMARK 2.2. One can derive similar inequalities as (2.4) when $\lambda_G$ is replaced with a higher eigenvalue $\lambda$ of $-L_G$ [cf. (1.7)] and $W$ in (2.3) with some connected component $U$ of the set $\{\psi > 0\} \subseteq G$, for some normalized eigenfunction $\psi$ of $-L_G$ attached to $\lambda$. Together with the invariance principle for the simple random walk on $\mathbb{Z}$, this yields quantitative lower bounds on the probability that the walk $X.$ travels in a cylinder $U \times \mathbb{Z}$ within time of order $\lambda^{-1}$ to a distance of order $\lambda^{-1/2}$ in the vertical direction, when starting at $x$ such that $g = \pi_G(x)$ corresponds to a value $\psi(g)$ "comparable" to $\max_U \psi$. In this fashion one obtains certain "escape routes" for the walk in the discrete cylinder $E$. In a way, the localization procedure we employ in the derivation of lower bounds on the disconnection time enables us to construct "easy escape routes" for the walk that only needs to travel in the vertical direction at distances of order $h'$ instead of distances of order $h$; see (2.12) and (3.4) below. It also avoids the use of detailed knowledge of the structure of higher eigenfunctions of $-L_G$.

We turn to the second result of this section that will be instrumental for the localization procedure. We now wish to consider stopping times defined by (1.14) corresponding to two distinct values of the parameter $h$ and the choice $u = 0$ (for simplicity). We thus consider [compare with (1.17)]

$$(2.12) \qquad 1 \leq h' \leq h = 2([\sqrt{t_G}(\log|G|)^2] + 1),$$

denote by $B', \widetilde{B}'$ and $B, \widetilde{B}$ the corresponding boxes when $u = 0$ [cf. (1.13)] as well as by $R'_k, D'_k, k \geq 1$, and $R_k, D_k, k \geq 1$, the corresponding stopping times; see (1.14). We also consider a subset of $G$, where the localization will take place:

$$(2.13) \qquad A \subseteq G.$$

We introduce the variables counting the visits of $X_{R'_k}, k \geq 1$, to $A \times \mathbb{Z}$ during the various intervals $[R_\ell, D_\ell - 1], \ell \geq 1$. We recall that in view of (2.12) all $R'_k$ occur during some $[R_\ell, D_\ell - 1], \ell \geq 1$. We thus define

$$U_1^A = \sum_{k \geq 1} 1\{X_{R'_k} \in A \times \mathbb{Z}, R'_k < D_1\} \circ \theta_{R_1}, \quad \text{and for } \ell \geq 1,$$

$$(2.14) \qquad U_\ell^A = U_1^A \circ \theta_{R_\ell}.$$

Clearly, the expectation under $P_x$ of $U_1^G$ only depends on $|\pi_\mathbb{Z}(x)|$, and we introduce

$$(2.15) \qquad \eta = E_x[U_1^G] \qquad \text{for } x \in G \times \{-h', h'\} \text{ arbitrary.}$$

Considering successive displacements at distance $h'$ of the simple random walk on $\mathbb{Z}$, that is, the iterates $\gamma_k$, $k \geq 0$, of the stopping time $\gamma = \inf\{n \geq 0, |Z_n - Z_0| = h'\}$,

$$\gamma_0 = 0, \qquad \gamma_{k+1} = \gamma \circ \theta_{\gamma_k} + \gamma_k \qquad \text{for } k \geq 0,$$

one knows that $\widehat{Z}_k = \frac{1}{h'} Z_{\gamma_k}$, $k \geq 0$, under $P_{h'}^\mathbb{Z}$ has the distribution of a simple random walk on $\mathbb{Z}$ starting at 1, that is, $P_1^\mathbb{Z}$. Using this identity, we see that $\eta$ is bounded from below by the expected number of successive returns to the interval $[-1, 1]$ of the one-dimensional simple random walk starting at 1, up to the exit time from $[-[\frac{h}{h'}], [\frac{h}{h'}]]$. Similarly, it is bounded from above by the expected number of returns to $[-1, 1]$ of a simple random walk starting at 1 up to the exit time from $[-[\frac{2h}{h'}] - 1, [\frac{2h}{h'}] + 1]$. With standard estimates on a simple random walk, it is straightforward to infer that, for some $c > 1$,

$$(2.16) \qquad \frac{1}{c} \frac{h}{h'} \leq \eta \leq c \frac{h}{h'}.$$

Finally, for $u \in \mathbb{Z}$, we denote with $\nu_u$ the equidistribution at level $u$ in $E$:

$$(2.17) \qquad \nu_u = \frac{1}{|G|} \sum_{x \in E, \pi_\mathbb{Z}(x) = u} \delta_x.$$

The second result of this section is (see above Theorem 1.2 for the terminology) the following:

PROPOSITION 2.3. *For large $G$, when $2 \leq \ell \leq |G|^2$, $v > 0$, $A \subseteq G$, and $x \notin \widetilde{B}$, one has the following:*

$$(2.18) \qquad P_x\left[U_1^A + \cdots + U_\ell^A > \ell \eta \frac{|A|}{|G|}(1 + v)\right] \leq c \exp\left\{-c \frac{|A|}{|G|} v(v \wedge 1)\ell\right\},$$

$$P_x\left[U_1^G > \frac{h}{h'} v\right] \leq 2 \exp\left\{-\frac{v}{2}\right\}$$

(2.19)

*for any $x \in E$, $v > 0$.*

PROOF. We begin with the proof of (2.18). We consider the variables $\overline{U}_\ell^A$, $\ell \geq 1$, attached to continuous time random walk $\overline{X}.$, obtained by replacing $R_k, D_k, R'_k, D'_k$ with $\overline{R}_k, \overline{D}_k, \overline{R}'_1, \overline{D}'_1$ in (2.14). The discrete skeleton of $\overline{X}.$ has the same law as $X.$ and therefore, under $P_x$, for arbitrary $x$ in $E$, $\overline{U}_\ell^A$, $\ell \geq 1$, has the same law as $U_\ell^A, \ell \geq 1$. As a result, for $x \notin \widetilde{B}$, $1 \leq \ell \leq |G|^2, A \subseteq G$, and $\lambda > 0$, the strong Markov property together with the above remark yields

$$E_x[\exp\{\lambda(U_1^A + \cdots + U_\ell^A)\}]$$
(2.20)
$$= E_x[\exp\{\lambda(\overline{U}_1^A + \cdots + \overline{U}_\ell^A)\}]$$
$$= E_x[\exp\{\lambda(\overline{U}_1^A + \cdots + \overline{U}_{\ell-1}^A)\}E_{\overline{X}_{\overline{D}_{\ell-1}}}[E_{\overline{X}_{\overline{R}_1}}[\exp\{\lambda \overline{U}_1^A\}]]],$$

where in the case $\ell = 1$, we use the convention $\overline{D}_0 = 0$, and the term before the inner expectation is omitted. We will use the following:

LEMMA 2.4. *For large $G$, where $x \notin \widetilde{B}$ and $z \in G \times \{-h, h\}$,*
(2.21)
$$P_x[\overline{X}_{\overline{R}_1} = z] \leq |G|^{-1}(1 + |G|^{-2})$$

*[of course the left-hand side vanishes unless $\pi_{\mathbb{Z}}(x)\pi_{\mathbb{Z}}(z) > 0$].*

PROOF. Without loss of generality, we assume that $x = (g, u), z = (f, h)$, with $u \geq 2h$, and $g, f \in G$. Using the exponential martingales $\exp\{\nu \overline{Z}_t - 2t \cosh \nu - 1)\}, t \geq 0$ (see, e.g., Lemma 3.2, page 175 of [9]), and applying Doob's inequality (see (2.46), page 63 of [9]), after optimization over $\nu \geq 0$, one obtains the (classical) bound

(2.22) $\quad P_0^{\mathbb{Z}}[\overline{H}_v \leq t] \leq \exp\left\{-cv \log\left(1 + c\frac{v}{t}\right)\right\} \qquad \text{for } v \in \mathbb{Z}_+, t > 0.$

With (1.5) and (1.11), we see that, for $t \geq t_G$

$$P_x[\overline{X}_{\overline{R}_1} = z] \quad = \quad P_g^G \otimes P_u^{\mathbb{Z}}[\overline{Y}_{\overline{H}_h^{\mathbb{Z}}} = z]$$
(2.23)
$$\leq \quad P_u^{\mathbb{Z}}[\overline{H}_h \leq t] + P_g^G \otimes P_u^{\mathbb{Z}}[\overline{Y}_{\overline{H}_h^{\mathbb{Z}}} = z, \overline{H}_h^{\mathbb{Z}} \geq t]$$
$$\overset{(2.22)}{\leq} \quad \exp\left\{-ch \log\left(1 + c\frac{h}{t}\right)\right\}$$
$$+ \int_t^\infty \frac{1}{|G|}\left(1 + \frac{1}{2}\exp\{-\lambda_G(s - t_G)\}\right)P_u^{\mathbb{Z}}[H_h \in ds].$$

We choose $t = 3t_G$. Observe that when $h/t \leq 1$, $\exp\{-ch\log(1+ch/t)\} \leq \exp\{-ch^2/t\} \leq \exp\{-c(\log|G|)^2\}$, with (2.12), and otherwise if $h/t \geq 1$, $\exp\{-ch\log(1+\frac{ch}{t})\} \leq \exp\{-ch\} \leq \exp\{-c(\log|G|)^2\}$. In addition, $\exp\{-2 \times \lambda_G t_G\} = 4|G|^{-2}$, due to (1.10). Coming back to (2.23), the claim (2.21) follows. □

As a result of Lemma 2.4 and the symmetry between positive and negative heights, for large $G$, the inner expectation in the last line of (2.20) is smaller than [in the notation of (2.17)]

(2.24) $\quad (1 + |G|^{-2})E_{\nu_h}[\exp\{\lambda \overline{U}_1^A\}] \leq (1 + |G|^{-2})E_{\nu_{h'}}[\exp\{\lambda \overline{U}_1^A\}],$

using the strong Markov property at time $\overline{R}'_1$ in the last step, and the fact that $\overline{X}_{\overline{R}'_1}$ is distributed as $\nu_{h'}$ under $P_{\nu_h}$. Iterating we see that for large $G$, for $x \notin \widetilde{B}$, $1 \leq \ell \leq |G|^2$, $A \subseteq G$ and $\lambda > 0$,

$$
\begin{aligned}
E_x[\exp\{\lambda(U_1^A + \cdots + U_\ell^A)\}] &\leq (1 + |G|^{-2})^\ell E_{\nu_{h'}}[\exp\{\lambda \overline{U}_1^A\}]^\ell \\
&\leq e E_{\nu_{h'}}[\exp\{\lambda \overline{U}_1^A\}]^\ell.
\end{aligned}
\tag{2.25}
$$

Observe also that using Taylor's formula with integral remainder to give a development to first order of the function $u \to e^{\lambda u}$, we find

$$
\begin{aligned}
E_{\nu_{h'}}[\exp\{\lambda \overline{U}_1^A\}] &= 1 + \lambda E_{\nu_{h'}}[\overline{U}_1^A] \\
&\quad + \lambda^2 E_{\nu_{h'}}\left[\int_0^1 ds \int_0^s dt (\overline{U}_1^A)^2 \exp\{\lambda t \overline{U}_1^A\}\right] \\
&\leq 1 + \lambda E_{\nu_{h'}}[\overline{U}_1^A] + \frac{\lambda^2}{2} E_{\nu_{h'}}[(\overline{U}_1^A)^2 \exp\{\lambda \overline{U}_1^A\}].
\end{aligned}
\tag{2.26}
$$

Since $\overline{\mu}$ is the stationary distribution of $\overline{Y}$, in view of (1.5) and (2.17) we find

$$
\begin{aligned}
E_{\nu_{h'}}[\overline{U}_1^A] &= E_{\overline{\mu}}^G \otimes E_{h'}^{\mathbb{Z}}\left[\sum_{k \geq 1} 1\{\overline{Y}_{\overline{R}'_k} \in A, \overline{R}'_k < \overline{D}_1\}\right] \\
&= E_{h'}^{\mathbb{Z}}\left[\sum_{k \geq 1} 1\{\overline{R}'_k < \overline{D}_1\}\right] \overline{\mu}(A) \stackrel{(2.14),(2.15)}{=} \eta \overline{\mu}(A),
\end{aligned}
\tag{2.27}
$$

with an abuse of notation when viewing the $\overline{R}'_k$, $k \geq 1$ and $\overline{D}_1$ as defined in terms of $\overline{Z}.$ alone. With analogous arguments, we also have

$$
\begin{aligned}
E_{\nu_{h'}}&[(\overline{U}_1^A)^2 \exp\{\lambda \overline{U}_1^A\}] \\
&\leq E_{\nu_{h'}}[(\overline{U}_1^A \overline{U}_1^G \exp\{\lambda \overline{U}_1^G\}] \\
&= E_{\overline{\mu}}^G \otimes E_{h'}^{\mathbb{Z}}\left[\sum_{k \geq 1} 1\{\overline{Y}_{\overline{R}'_k} \in A, \overline{R}'_k < \overline{D}_1\}, \overline{U}_1^G \exp\{\lambda \overline{U}_1^G\}\right] \\
&= \overline{\mu}(A) E_{h'}^{\mathbb{Z}}\left[\sum_{k \geq 1} 1\{\overline{R}'_k < \overline{D}_1\} \overline{U}_1^G \exp\{\lambda \overline{U}_1^G\}\right] \\
&= \overline{\mu}(A) E_{h'}^{\mathbb{Z}}[(\overline{U}_1^G)^2 \exp\{\lambda \overline{U}_1^G\}].
\end{aligned}
\tag{2.28}
$$

Further, note that the simple random walk on $\mathbb{Z}$ starting at $2h'$ reaches $2h$ before $h'$ with probability $h'/(2h - h') \geq h'/(2h)$. With a repeated application of the strong Markov property, we find that

$$
P_x[U_1^G \geq m] \leq P_x[R'_m < D_1] \leq \left(1 - \frac{1}{2}\frac{h'}{h}\right)^{m-1} \quad \text{for } x \in E, m \geq 1.
\tag{2.29}
$$

Hence, for a suitable small enough positive constant $c_3$, and any $x \in E$,

$$(2.30) \qquad E_x\left[\exp\left\{c_3 \frac{h'}{h} U_1^G\right\}\right] = E_x\left[\exp\left\{c_3 \frac{h'}{h} \overline{U}_1^G\right\}\right] \le 2.$$

Combining (2.26), (2.27), (2.28) and (2.30), we see that for $\lambda \le \frac{c_3}{2} \frac{h'}{h}$ one has

$$(2.31) \qquad E_{\nu_{h'}}[\exp\{\lambda \overline{U}_1^A\}] \le 1 + \eta \overline{\mu}(A)\lambda + c\left(\frac{h}{h'}\right)^2 \overline{\mu}(A)\lambda^2.$$

Returning to (2.25), we obtain for $v > 0$, and $\lambda \le \frac{c_3}{2} \frac{h'}{h}$,

$$P_x[U_1^A + \cdots + U_\ell^A > \eta \overline{\mu}(A)(1+v)\ell]$$
$$\le \exp\left\{-\lambda \eta \overline{\mu}(A)(1+v)\ell + 1 + \ell\left(\lambda \eta \overline{\mu}(A) + c\left(\frac{h}{h'}\right)^2 \overline{\mu}(A)\lambda^2\right)\right\}$$
$$\stackrel{(2.16)}{\le} \exp\left\{1 - \ell \overline{\mu}(A)\left[c \frac{\lambda h}{h'} v - c'\left(\frac{\lambda h}{h'}\right)^2\right]\right\}.$$

Optimizing over $\lambda$ with the definition (1.6), we obtain (2.18). We now turn to the proof of (2.19). With (2.29), we find that for any $x \in E$, $v > 0$,

$$P_x\left[U_1^G > \frac{h}{h'} v\right] \le \left(1 - \frac{1}{2}\frac{h'}{h}\right)^{[(h/h')v]} \le 2\exp\left\{-\frac{1}{2}\frac{h'}{h}\frac{h}{h'} v\right\} = 2\exp\left\{-\frac{v}{2}\right\},$$

whence (2.19).  □

We then continue with the third result of this section. It provides bounds that will be especially helpful when $h$ in (2.12) is not too large [i.e., $\lambda_G$ large enough; see (3.34) and (4.10)]. We use the terminology introduced above Theorem 1.2 for the next lemma.

LEMMA 2.5. *For large $G$,*

$$(2.32) \quad P_x[H_{x'} < T_{\widetilde{B}}] \le c\frac{h}{|G|} \qquad \text{for } x \in G \times \{-h, h\} \quad \text{and} \quad |\pi_{\mathbb{Z}}(x')| \le \frac{h}{2}.$$

*Moreover, for any $V \subseteq G$ with $|V|h \le |G|(\log|G|)^{-2}$ and $x \in E$, one has the following:*

$$(2.33) \quad \begin{aligned} &\text{(i)} \quad E_x\left[\exp\left\{\frac{c}{\sqrt{t_G}}|\pi_{\mathbb{Z}}(X_{[0,T_{\widetilde{B}}-1]} \cap (V \times \mathbb{Z}))|\right\}\right] \le 2, \\ &\text{(ii)} \quad E_x\left[\exp\left\{\frac{c}{t_G}|\pi_G(X_{[0,T_{\widetilde{B}}-1]} \cap (V \times \mathbb{Z}))|\right\}\right] \le 2. \end{aligned}$$

PROOF. We begin with the proof of (2.32). The argument resembles the proof of (2.21). Without loss of generality, we assume that $x = (g, h)$, $x' = (g', u')$, with $|u'| \leq \frac{h}{2}$, and write, using similar bounds as in (2.22) and below (2.23),

$$P_x[H_{x'} < T_{\widetilde{B}}] = P_x[\overline{H}_{x'} < \overline{T}_{\widetilde{B}}]$$

$$\leq P_x[\overline{H}_{x'} \leq t_G] + P_x[t_G < \overline{H}_{x'} < \overline{T}_{\widetilde{B}}]$$

$$\leq \exp\{-c(\log |G|)^2\} + P_x[t_G < \overline{H}_{x'} < \overline{T}_{\widetilde{B}}]$$

$$\overset{(1.5)}{\leq} \exp\{-c(\log |G|)^2\}$$

(2.34)
$$+ \int_{t_G}^{\infty} c P_g^G[\overline{Y}_t = g'] P_h^{\mathbb{Z}}[\overline{Z}_t = u', \overline{T}_{\widetilde{I}} > t] dt$$

$$\overset{(1.11)}{\leq} \exp\{-c(\log |G|)^2\} + \frac{c}{|G|} \int_{t_G}^{\infty} P_h^{\mathbb{Z}}[\overline{Z}_t = u', \overline{T}_{\widetilde{I}} > t] dt$$

$$\leq \exp\{-c(\log |G|)^2\} + \frac{c}{|G|} E_0^{\mathbb{Z}} \left[ \int_0^{\overline{T}_{[-4h, 4h]}} 1_{\{\overline{Z}_t = 0\}} dt \right]$$

$$\leq c \frac{h}{|G|},$$

whence (2.32). As for (2.33)(i), we first note that

(2.35) $\quad |\pi_{\mathbb{Z}}(X_{[0, T_{\widetilde{B}} - 1]} \cap (V \times \mathbb{Z}))| = \sum_{|\widetilde{u}| < 2h} 1\{H_{V \times \{\widetilde{u}\}} < T_{\widetilde{B}}\}.$

With an argument similar to Khasminskii's lemma [12] (see also, e.g., (2.46) of [8]) the claim (2.33) (i) follows once we show that

(2.36) $\quad \sup_{x \in E} E_x \left[ \sum_{|\widetilde{u}| < 2h} 1\{H_{V \times \{\widetilde{u}\}} < T_{\widetilde{B}}\} \right] \leq c \sqrt{t_G}.$

To prove (2.36), note that the above expectation is equal to

$$\sum_{|\widetilde{u}| < 2h} P_x[\overline{H}_{V \times \{\widetilde{u}\}} < \overline{T}_{\widetilde{B}}] \leq \sum_{|\widetilde{u}| < 2h} P_x[\overline{H}_{V \times \{\widetilde{u}\}} < t_G]$$

$$+ \sum_{|\widetilde{u}| < 2h} P_x[t_G < \overline{H}_{V \times \{\widetilde{u}\}} < \overline{T}_{\widetilde{B}}]$$

$$\leq 2 E_u^{\mathbb{Z}} \left[ \sup_{0 \leq s \leq t_G} |\overline{Z}_s - \overline{Z}_0| + 1 \right]$$

$$
\begin{aligned}
(2.37) \qquad & + cE_x\left[\int_{t_G}^{\infty} 1\{\overline{X}_t \in V \times \mathbb{Z}, t < T_{\widetilde{B}}\}\,dt\right] \\
& \stackrel{(1.5)}{\leq} c\sqrt{t_G} + c\int_{t_G}^{\infty} P_g^G[\overline{Y}_t \in V]P_u^{\mathbb{Z}}[t < \overline{T}_{\widetilde{I}}]\,dt \\
& \stackrel{(1.11)}{\leq} c\sqrt{t_G} + c\frac{|V|}{|G|}E_u^{\mathbb{Z}}[\overline{T}_{\widetilde{I}}] \\
& \leq c\sqrt{t_G} + c|V|\frac{h^2}{|G|}.
\end{aligned}
$$

However, with (2.12) and our assumption on $V$,

$$\frac{|V|h^2}{|G|} \leq c(\log|G|)^2 \frac{|V|h}{|G|}\sqrt{t_G} \leq c\sqrt{t_G},$$

and (2.36) follows.

The proof of (2.33)(ii) is similar. In place of (2.36), we have to check that

$$(2.38) \qquad \sup_{x \in E} E_x\left[\sum_{\widetilde{g} \in V} 1\{H_{\{\widetilde{g}\} \times \widetilde{I}} < T_{\widetilde{B}}\}\right] \leq ct_G.$$

Moreover, since

$$\sum_{\widetilde{g} \in V} E_x[t_G < \overline{H}_{\{\widetilde{g}\} \times \widetilde{I}} < \overline{T}_{\widetilde{B}}] \leq cE_x\left[\int_{t_G}^{\infty} 1\{\overline{X}_t \in V \times \mathbb{Z}, t < \overline{T}_{\widetilde{B}}\}\,dt\right],$$

the claim (2.38) follows from a straightforward modification of (2.37). This concludes the proof of (2.33). □

We now turn to the fourth and last result of this section, which highlights the strategy we will employ when bounding the disconnection time from below. We depart from the line of attack in [8], which was based on the fact that a finite subset $S$ disconnecting $(\mathbb{Z}/N\mathbb{Z})^d \times \mathbb{Z}$ must somewhere be "locally big," thanks to isoperimetric controls; compare Lemmas 2.4 and 2.5 of [8]. Here instead we construct paths that prevent disconnection.

To this end, we consider integers $M, M', L' \geq 1$, as well as

(2.39) $\qquad V_i, 1 \leq i \leq M \qquad$ nonempty connected subsets of $G$,

$$\text{with } \bigcup_{i=1}^{M} V_i \text{ connected.}$$

We define for $1 \leq j \leq M'$ the intervals of $\mathbb{Z}$:

$$(2.40) \quad J_j = [(j-1)L', jL'], \qquad 1 \leq j \leq M' \qquad (\text{so } |J_j| = L'+1),$$

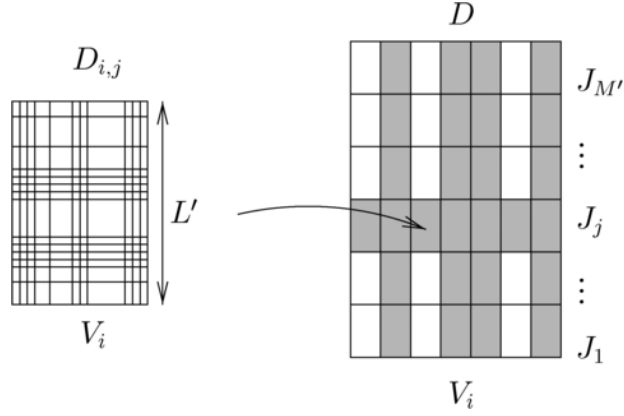

FIG. 1. *A schematic illustration of "S is good for D" and "S is thin in $D_{ij}$." The sets corresponding to the thin lines in the box on the left are in the complement of S. Also S is thin in the shaded boxes on the right-hand side of the figure.*

as well as the subsets of $E$:

$$
\begin{aligned}
D_{i,j}(u) &= V_i \times (J_j + u), \\
D(u) &= \bigcup_{1 \leq i \leq M, 1 \leq j \leq M'} D_{i,j}(u) \quad \text{for } u \in \mathbb{Z}.
\end{aligned}
\tag{2.41}
$$

We simply write $D_{i,j}$ and $D$ when $u = 0$.

Given a finite subset $S \subset E$, we say that *S is thin in $D_{i,j}(u)$*, when

$$(2.42) \quad |\pi_G(D_{i,j}(u) \cap S)| < \frac{|V_i|}{2} \quad \text{and} \quad |\pi_\mathbb{Z}(D_{i,j}(u) \cap S)| < \frac{|J_j|}{2}\left(= \frac{L'+1}{2}\right).$$

Using a type of renormalized version of (2.42), we say that *S is good for $D(u)$*, when

$$
\begin{aligned}
\text{(i)} \quad &|\{i \in [1, M]; S \text{ is not thin in } D_{i,j}(u) \text{ for some } j \in [1, M']\}| < \frac{M}{2}, \\
\text{(ii)} \quad &|\{j \in [1, M']; S \text{ is not thin in } D_{i,j}(u) \text{ for some } i \in [1, M]\}| < M'.
\end{aligned}
\tag{2.43}
$$

We formulated (2.43) in a way which highlights the analogy with (2.42); in particular, "S not thin in $D_{i,j}(u)$" for $(i,j) \in [1,M] \times [1,M']$ is the counterpart of "x in S" for $x \in D_{i,j}(u)$ in (2.42). Note that when $S \cap D(u) = \varnothing$, $S$ is automatically good for $D(u)$ (see Figure 1). Our last result in this section is the following:

PROPOSITION 2.6. $(n \geq 0)$

$$(2.44) \quad \{\mathcal{T}_G \leq n\} \subseteq \{ \text{for some } u \in \mathbb{Z}; X_{[0,n]} \text{ is not good in } D(u) \}.$$

PROOF. We prove (2.44) by contradiction. We denote with $\mathcal{G}_n$ the complement of the event in the right-hand side of (2.44). We fix a trajectory in $\mathcal{G}_n$, and set $S = X_{[0,n]}$. We choose $k_- \leq k_+$ in $\mathbb{Z}$, so that

$$(2.45) \quad \left(\bigcup_{i=1}^M V_i\right) \times [\min \pi_{\mathbb{Z}}(S), \max \pi_{\mathbb{Z}}(S)] \subseteq \bigcup_{k_- \leq k \leq k_+} D(kM'L').$$

We say that $D_{i,j}(kM'L')$ belongs to a thin column, respectively to a thin row, of $D(kM'L')$ when $S$ is thin in each $D_{i,j'}(kM'L')$, $1 \leq j' \leq M'$, respectively each $D_{i',j}(kM'L')$, $1 \leq i' \leq M$. The first observation is the following:

(2.46) Any two boxes $D_{i_-,1}(k_-M'L')$ and $D_{i_+,M'}(k_+M'L')$ in a thin column of $D(k_-M'L')$ and $D(k_+M'L')$ respectively can be linked by a path of boxes $\mathcal{D}_\ell = D_{i_\ell,j_\ell}(k_\ell M'L')$, $0 \leq \ell \leq m$,

such that:

(i) the path starts in $D_{i_-,1}(k_-M'L')$ and ends in $D_{i_+,M'}(k_+M'L')$,
(ii) for each $0 \leq \ell \leq m$, $S$ is thin in $\mathcal{D}_\ell$,
(iii) for each $0 \leq \ell < m$, either the boxes $\mathcal{D}_\ell$ and $\mathcal{D}_{\ell+1}$ are vertically abutting, that is, $i_\ell = i_{\ell+1}$ and $|j_\ell L' + k_\ell M'L' - j_{\ell+1}L' - k_{\ell+1}M'L'| = L'$, or side-wise abutting, that is, $j_\ell = j_{\ell+1}, k_\ell = k_{\ell+1}$, and $V_{i_\ell} \cap V_{i_{\ell+1}} \neq \varnothing$.

Indeed, since the trajectory of the walk belongs to $\mathcal{G}_n$, any two boxes in thin columns of $D(kM'L')$ can be linked by such a nearest neighbor path of boxes in a thin row or thin column of $D(kM'L')$; see (2.43)(ii). In addition, with (2.43)(i), for any $k$, at least one $1 \leq i \leq M$ is such that $D_{i,L'}(kM'L')$ and $D_{i,1}((k+1)M'L')$ are both in thin columns of $D(kM'L')$ and $D((k+1)M'L')$, respectively. The claim (2.46) follows.

The next observation is that for $0 \leq \ell < m$:

(2.47) any $z, z'$ respectively in $\mathcal{D}_\ell$ and $\mathcal{D}_{\ell+1}$, such that $\pi_{\mathbb{Z}}(z) \notin \pi_{\mathbb{Z}}(S \cap \mathcal{D}_\ell)$ or $\pi_G(z) \notin \pi_G(S \cap \mathcal{D}_\ell)$ and a similar condition for $z'$ with $\ell$ replaced by $\ell + 1$, can be joined by a nearest neighbor path in $(\mathcal{D}_\ell \cup \mathcal{D}_{\ell+1}) \setminus S$.

[Such points exist in view of (2.46)(ii) and (2.42).]

Indeed, one can construct a path within $\mathcal{D}_\ell \setminus S$ or $\mathcal{D}_{\ell+1} \setminus S$ between two points in the same box $\mathcal{D}_\ell$ or $\mathcal{D}_{\ell+1}$ that satisfies the above mentioned property thanks to (2.42) and the fact that each $V_i, 1 \leq i \leq M$, is connected. Then using the fact that $\mathcal{D}_\ell$ and $\mathcal{D}_{\ell+1}$ are either vertically abutting with $V_{i_\ell} = V_{i_{\ell+1}}$, or sidewise abutting with $V_{i_\ell} \cap V_{i_{\ell+1}} \neq \varnothing$, in view of (2.42), we can either find $g \in V_{i_\ell} = V_{i_{\ell+1}}$ such that $g \in \pi_G((\mathcal{D}_\ell \cup \mathcal{D}_{\ell+1}) \setminus S)$, or $u \in \pi_{\mathbb{Z}}((\mathcal{D}_\ell \cup \mathcal{D}_{\ell+1}) \setminus S)$. The claim (2.47) readily follows in view of the previous remark.

With (2.46) and (2.47), we can then find a nearest neighbor path in $E \backslash S$ which starts at a point having $\mathbb{Z}$-projection equal to $\min \pi_{\mathbb{Z}}(S)$ and ends at a point having $\mathbb{Z}$-projection equal to $\max \pi_{\mathbb{Z}}(S)$. Hence, $S$ does not disconnect $E$. In other words, we have shown that $\mathcal{G}_n \subseteq \{\mathcal{T}_G > n\}$, and (2.44) follows. □

In Section 4, when working in the presence of a transient pocket (cf. Theorem 4.1), we use the simple case $M = M' = 1$, so that "$S$ good for $D(u)$" means that $S$ is thin in $D_{1,1}(u)$. On the other hand, in order to handle the possible presence of recurrent pockets in Section 5, we use the above Proposition 2.6 with $M$ and $M' > 1$; see Theorem 5.2.

**3. Localization technique and rarefied excursions.** We develop the localization technique in this section. We focus on what happens in a certain "pocket" $A$ of $G$ (see Figure 2), where we have control over the decay of the killed heat kernel; see (3.4). We are interested in the excursions performed before time $|G|^{2(1-\delta)}$, corresponding to successive entrances of the walk in a not too big box $C$ with $G$ projection denoted by $V$, "well inside $A$," and departures of the walk at distances of order $h' \leq h$ in the $\mathbb{Z}$-direction. We pick $h'$ [cf (3.6)] so that, on the one hand, it is large enough and thereby makes it rare to hit $C$ when starting at vertical distance of order $h'$ from $C$, and on the other hand, small enough so that in the later contexts of Sections 4 and 5, we are able to check (3.5), and thereby discard what happens outside $A \times \mathbb{Z}$, when analyzing these excursions. Our key result (cf. Proposition 3.2) shows that for our purpose we can assume that only a finite number of excursions take place. This is instrumental when later constructing connections with Proposition 2.6. Our convention on constants for this section appears above Remark 3.1. We first introduce some definitions.

Recalling $d_0$ from (1.2), we denote with $\mathcal{G}_0(d_0, \delta, \delta', \gamma, a)$, where $0 < \delta < \frac{1}{2}$, $0 < \delta' < \frac{\delta}{8}$, $\gamma > \frac{1}{2}$, $a > 0$, the class of finite connected graphs $G$ satisfying (1.2), such that either

$$\lambda_G \leq |G|^{-2(1-\delta)}(\log |G|)^{-1}, \tag{3.1}$$

or (3.1) does not hold and there exist

$$G \supseteq A \supseteq V \quad \text{with } V \text{ connected, so that} \tag{3.2}$$

$$|A| \geq \lambda_G^{-1/2} |G|^{15\delta/16}, \qquad |G|^{\delta/8} \geq |V| \geq \frac{1}{d_0} |G|^{\delta'} \tag{3.3}$$

and

$$P_g^G[Y_n = g, n < T_A] \leq \frac{a}{n^\gamma}, \qquad \text{for } g \in A, n \geq 1. \tag{3.4}$$

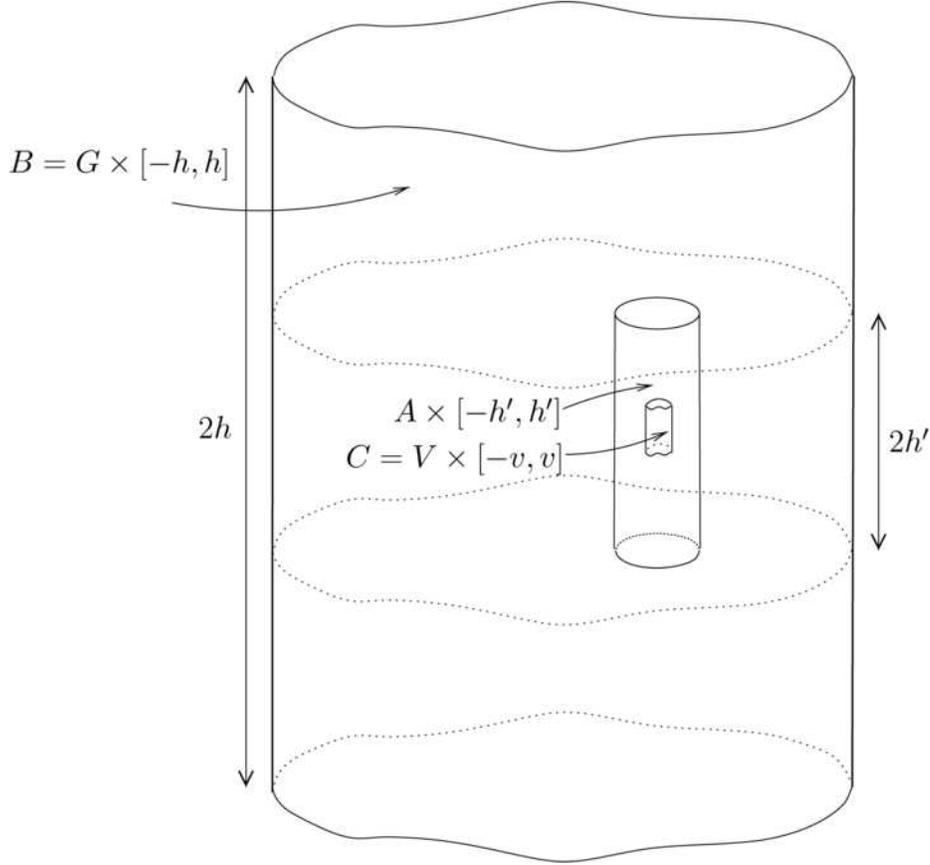

Fig. 2. *A drawing of some of the sets that appear in (3.2)–(3.7).*

In addition, we assume that

$$\text{(3.5)} \qquad \frac{|G|^{1-\delta/8}}{h'} \sup_{x \in A^c \times \mathbb{Z}} P_x[H_C < T_{\widetilde{B}'}] \leq 1,$$

where we have set, with $h$ as in (2.12),

$$\text{(3.6)} \qquad h' = [(\sqrt{|A|}|G|^{-\delta/8})^{1/\gamma}] \wedge h$$

and

$$\text{(3.7)} \quad C(u) = V \times (u + [-v, v]) \qquad \text{with } v = |V| \wedge \left[\frac{h'}{2}\right] \text{ and } u \in \mathbb{Z},$$

[and we write $C$ in place of $C(0)$]. We also consider the classes

$$\text{(3.8)} \qquad \mathcal{G}_1(d_0, \delta, \delta', \gamma, a) = \{G \in \mathcal{G}_0; (3.1) \text{ does not hold}\},$$

$$(3.9) \quad \mathcal{G}_{\text{ext}}(d_0, \delta) = \{G \text{ satisfying } (1.2) \text{ such that either } (3.1)$$
$$\text{holds or } \lambda_G \geq |G|^{-\delta/10}\}.$$

This latter definition corresponds to "extreme values of $\lambda_G$" [cf. (1.9)] and this class will be considered in Theorem 4.3. Throughout the section we use similar conventions concerning positive constants or the expression "for large $G$," as mentioned below (1.14), except that constants may now depend on $d_0, \delta, \delta', \gamma, a$ (and not just $d_0$). Let us give some comments about these parameters. The parameter $\delta$ ultimately measures the quality of the lower bound we derive on the disconnection time, with a similar interpretation as in (0.9). The numbers $\gamma$ and $a$ control the on-diagonal decay of the transition kernel of the walk on $G$ killed when exiting the "pocket" $A$, whereas $\delta'$ ensures the nondegeneracy of $V$ the $G$-projection of the box $C$ sitting inside $A \times [-h', h']$. The choice of $h'$ in (3.6), as well as (3.5), addresses the two conflicting constraints expressed at the beginning of this section. With the first and second inequality of (3.3), $A$ is not too small in $G$, $V$ is small in $A$, and due to (3.5), sits "well inside" $A$. The last inequality of (3.3) enforces a lower bound on $|V|$ which will ensure that we are not looking at too small a scale in $G$, and the multiplicity of boxes $C(u), |u| \leq |G|^2$, we later need to consider does not beat the probabilistic estimates we derive; see, for instance, (4.9), above and below (4.10)(ii), as well as the last line of (5.25). Let us mention that in some applications the values of $\gamma$ and $a$ will be fixed (see, e.g., Corollary 4.6 and Corollary 5.3) but in Corollary 4.5 we let $\gamma$ and $a$ depend on $\delta$ [and tend to infinity as $\delta$ tends to 0; cf. (4.20)].

REMARK 3.1. As a routine consequence of (1.2) and (3.4), one has

$$(3.10) \quad \text{for all } g, g' \in A, n \geq 1, \quad P_g^G[Y_n = g', n < T_A] \leq \frac{c}{n^\gamma}.$$

Indeed, $r_n(g, g') = P_g^G[Y_n = g', n < T_A]\mu(g')^{-1}$ is a symmetric function of $g, g'$, thanks to reversibility. Further, with the Chapman–Kolmogorov equations and Cauchy–Schwarz's inequality, one finds for $k \geq 1, g, g' \in A$,

$$r_{2k}(g, g') \leq r_{2k}(g, g)^{1/2} r_{2k}(g', g')^{1/2} \overset{(3.4)}{\leq} \frac{a}{(2k)^\gamma}(\mu(g)\mu(g'))^{-1/2},$$

$$r_{2k+1}(g, g') \leq r_{2k}(g, g)^{1/2} r_{2k+2}(g', g')^{1/2} \overset{(3.4)}{\leq} \frac{a}{(2k)^\gamma}(\mu(g)\mu(g'))^{-1/2},$$

whence (3.10).

We wish to control excursions consisting of successive returns to $C(u)$ and departures from $\widetilde{B}'(u) (= G \times \widetilde{I}'(u)$, with $\widetilde{I}'(u) = u + [-(2h'-1), 2h'-1])$.

To this end, we introduce for $u \in \mathbb{Z}$ the sequence of stopping times

$$\widetilde{R}_1^u = H_{C(u)}, \qquad \widetilde{D}_1^u = T_{\widetilde{B}'(u)} \circ \theta_{\widetilde{R}_1^u} + \widetilde{R}_1^u, \quad \text{and for } k \geq 1,$$

(3.11)

$$\widetilde{R}_{k+1}^u = \widetilde{R}_1^u \circ \theta_{\widetilde{D}_k^u} + \widetilde{D}_k^u, \qquad \widetilde{D}_{k+1}^u = \widetilde{D}_1^u \circ \theta_{\widetilde{D}_k^u} + \widetilde{D}_k^u,$$

and with a similar convention as below (1.14), we drop the superscript $u$ when $u = 0$, and simply write $\widetilde{R}_k, \widetilde{D}_k$. The next proposition shows the rarefaction of excursions between $C(u)$ and $E \setminus \widetilde{B}'(u)$, up to time $|G|^{2(1-\delta)}$, when $G$ is large in $\mathcal{G}_1$. It plays an important role in the present and next section.

PROPOSITION 3.2. *There is a positive constant $K_0$ (cf. above Remark 3.1), such that*

(3.12) $$\lim_{|G| \to \infty, G \in \mathcal{G}_1} \sup_{x \in E} P_x \left[ \sup_{u \in \mathbb{Z}} \sum_{k \geq 1} 1\{\widetilde{R}_k^u < |G|^{2(1-\delta)}\} > K_0 \right] = 0.$$

PROOF. We define [compare with (1.29); note that the value of $h$ in (1.29) is set by (1.17), whereas in the present section it is defined by (2.12)]

(3.13) $$m_* = \left[ \frac{|G|^{1-\delta}}{h} (\log |G|)^5 \right] + 1.$$

With (1.10), (2.12) and (3.8), we see that

(3.14)
$$\text{for large } G \text{ in } \mathcal{G}_1, \quad m_* \geq c(\log |G|)^2 \quad \text{and}$$
$$m_* \frac{|A|}{|G|} |G|^{\delta/16} \geq c(\log |G|)^{5/2}.$$

Recall the notation $R_k^u, D_k^u$ from below (1.14), with a choice for $h$ made in (2.12). With similar arguments as below (2.11) of [8], for large $G$ in $\mathcal{G}_1$ and arbitrary $z \in E$, $u \in \mathbb{Z}$, we find,

(3.15)
$$P_z[R_{m_*}^u \leq |G|^{2(1-\delta)}]$$
$$\leq P_0^{\mathbb{Z}}[H_{(m_*-1)h} \leq |G|^{2(1-\delta)}] \leq \left(1 - \frac{c}{|G|^{1-\delta}}\right)^{(m_*-1)h}$$
$$\leq \exp\{-c(\log |G|)^5\}.$$

We now turn to the control of the sum in the probability in (3.12) for $u = 0$, under $P_x$, with $x$ arbitrary in $E$. The case $u \neq 0$ will then follow using translation invariance in the $\mathbb{Z}$-direction. We first note that

$$(3.16) \qquad m_* \frac{h}{h'} \frac{|A|}{|G|} |G|^{\delta/16} \leq c \frac{|A|}{h'} |G|^{-15\delta/16} (\log |G|)^5 \stackrel{\text{def}}{=} k_*(A),$$

and also introduce the notation $k_*(G)$, when $G$ replaces $A$ in (3.16). With the help of Proposition 2.3, we can control the number of returns $R'_k$ with $X_{R'_k}$ in $A$ (or in $G$) that occur before $D_{m^*}$. Indeed, when $G$ is large in $\mathcal{G}_1$ and any $x \in E$, one has

$$P_x[U_1^A + \cdots + U_{m_*}^A > k_*(A)]$$

$$(3.17) \qquad \begin{aligned} &\leq\ P_x\Big[U_1^G > \frac{1}{2}k_*(A)\Big] + P_x\Big[U_2^A + \cdots + U_{m_*}^A > \frac{1}{2}k_*(A)\Big] \\ &\stackrel{\substack{(2.18),(2.19) \\ (2.16),(3.16)}}{\leq} 2\exp\Big\{-cm_* \frac{|A|}{|G|} |G|^{\delta/16}\Big\} + c\exp\Big\{-cm_* \frac{|A|}{|G|} |G|^{\delta/16}\Big\} \\ &\stackrel{(3.14)}{\leq}\ c\exp\{-c(\log|G|)^{5/2}\}. \end{aligned}$$

The same argument shows that

$$(3.18) \qquad \begin{aligned} P_x[U_1^G + \cdots + U_{m_*}^G > k_*(G)] &\leq\ c\exp\{-cm_*|G|^{\delta/16}\} \\ &\stackrel{(3.14)}{\leq}\ \exp\{-c|G|^{\delta/16}\}. \end{aligned}$$

With (3.15), (3.17) and (3.18), we have a bound on the number of returns $R'_k$ with $X_{R'_k}$ in $A$ (or in $G$) that occur before time $|G|^{2(1-\delta)}$. We then need to bound the probability of entering $C$ before exiting $\widetilde{B}'$. When $h' = h$ [cf. (3.6)], we will rely on (2.32) of Lemma 2.5. On the other hand when $h' < h$, we will use the following:

LEMMA 3.3 $[\gamma > \frac{1}{2}, \varepsilon \in (0, 2\gamma - 1)]$.  For $G$ in $\mathcal{G}_1$ with $|G| \geq c(\varepsilon)$, one has

$$(3.19) \quad P_x[H_{x'} < T_{(A \times \mathbb{Z}) \cap \widetilde{B}'}] \leq h'^{-(2\gamma - 1 - \varepsilon)} \qquad \text{for } x \in A \times \{-h', h'\}, x' \in C.$$

PROOF.  Pick $\delta_1 < \frac{1}{2}$. With no loss of generality, we assume that $x = (g, h')$, with $g \in A$ and $x' = (g', u') \in C$, so that $g' \in V$ and $|u'| \leq \frac{h'}{2}$; cf. (3.7). When $G$ is in $\mathcal{G}_1$, we thus find

$$(3.20) \qquad \begin{aligned} P_x[H_{x'} &< T_{(A \times \mathbb{Z}) \cap \widetilde{B}'}] \\ &=\ P_x[\overline{H}_{x'} < \overline{T}_{(A \times \mathbb{Z}) \cap \widetilde{B}'}] \end{aligned}$$

$$\stackrel{(1.5)}{\leq} P_{h'}^{\mathbb{Z}}[\overline{H}_{u'} < h'^{2(1-\delta_1)}]$$
$$+ cE_g^G \otimes E_{h'}^{\mathbb{Z}}\left[\int_{h'^{2(1-\delta_1)}}^{\infty} 1\{\overline{Y}_t = g', \overline{T}_A > t\}1\{\overline{Z}_t = u'\}\,dt\right].$$

Note that from (3.3), (3.6), (1.9) and (2.12), $\lim_{|G|\to\infty, G\in\mathcal{G}_1} h' = \infty$. With (2.22), for large $G$ in $\mathcal{G}_1$, the first term in the right-hand side of (3.20) is smaller than $\exp\{-ch'\log(1+ch'^{(1-2(1-\delta_1))})\} \leq \exp\{-ch'^{2\delta_1}\}$. We now turn to the last term of (3.20). We introduce the continuous piecewise linear increasing processes

$$A_t = \int_0^t \deg(\overline{Y}_r)\,dr, \qquad t \geq 0, \quad \text{and}$$

$$\tau_s = \inf\{t > 0; A_t > s\} = (A^{-1})_s, \qquad s \geq 0,$$

as well as the time changed process

$$\widetilde{Y}_s = \overline{Y}_{\tau_s}, \qquad s \geq 0.$$

Under $P_g^G$, $\widetilde{Y}.$ is the continuous walk on $G$ with the same discrete skeleton as $\overline{Y}.$, but with constant jump rate one. Note also that

$$\tau_s = \int_0^s \deg(\widetilde{Y}_r)^{-1}\,dr.$$

Performing the change of variable $t = \tau_s$, the last term of (3.20) equals

$$cE_g^G \otimes E_u^{\mathbb{Z}}\left[\int_{A_{h'^{2(1-\delta_1)}}}^{\infty} 1\{\widetilde{Y}_s = g', T_A^{\widetilde{Y}} > s\}1\{\overline{Z}_{\tau_s} = u'\}\deg(\widetilde{Y}_s)^{-1}\,ds\right]$$

$$\leq c\int_{h'^{2(1-\delta_1)}}^{\infty} E_g^G[\widetilde{Y}_s = g', T_A^{\widetilde{Y}} > s, q(\tau_s, u, u')]\,ds,$$

with $q(t, u, u') = P_u^{\mathbb{Z}}[\overline{Z}_t = u']$, and otherwise hopefully obvious notation. Note that when $N_t, t \geq 0$, is a Poisson counting process with rate 1,

(3.21) $$P\left[N_t < \frac{t}{2} \text{ or } N_t > 2t\right] \leq 2e^{-ct}, \qquad t > 0,$$

as follows from Cramér-type bounds. With (3.10) and a similar bound for $q(t, u, u')$, the last term of (3.20) is thus smaller than

$$c\int_{h'^{2(1-\delta_1)}}^{\infty} (t^{-(\gamma+1/2)} + e^{-ct})\,dt \leq c(h'^{(1-2\gamma)(1-\delta_1)} + e^{-ch'^{2(1-\delta_1)}}).$$

Choosing $\delta_1 < c(\varepsilon)$, so that $(2\gamma - 1)(1 - \delta_1) > 2\gamma - 1 - \varepsilon$, our claim (3.19) follows. $\square$

We now return to the task of bounding the sum in (3.12) for $u = 0$. We first analyze the case when [cf. (3.6)]

(3.22) $$h' < h.$$

For large $G$ in $\mathcal{G}_1$, under (3.22), for any $K \geq 2$ and $x \in E$,

$$P_x\left[\sum_{k \geq 1} 1\{\widetilde{R}_k \leq |G|^{2(1-\delta)}\} > K\right]$$

$$\leq P_x\left[\sum_{\ell \geq 1} 1\{H_C < T_{\widetilde{B}'}\} \circ \theta_{R'_\ell} 1\{R'_\ell \leq |G|^{2(1-\delta)}\} > K\right]$$

$$\stackrel{(3.15)}{\leq} \exp\{-c(\log|G|)^5\}$$

$$+ P_x\left[\sum_{\substack{1 \leq m \leq m_* \\ 1 \leq k}} 1\{H_C < T_{\widetilde{B}'}\} \circ \theta_{R'_k} 1\{R_m \leq R'_k < D_m\} > K\right]$$

(3.23)
$$\leq \exp\{-c(\log|G|)^5\}$$

$$+ P_x\Bigg[\sum_{\substack{1 \leq m \leq m_* \\ 1 \leq k}} 1\{H_C < T_{\widetilde{B}'}, \pi_G(X_0) \in A\} \circ \theta_{R'_k}$$

$$\times 1\{R_m \leq R'_k < D_m\} > \frac{K}{2}\Bigg]$$

$$+ P_x\Bigg[\sum_{\substack{1 \leq m \leq m_* \\ 1 \leq k}} 1\{H_C < T_{\widetilde{B}'}, \pi_G(X_0) \notin A\} \circ \theta_{R'_k}$$

$$\times 1\{R_m \leq R'_k < D_m\} > \frac{K}{2}\Bigg]$$

$$\stackrel{(3.17),(3.18)}{\leq} c\exp\{-c(\log|G|)^{5/2}\} + a_1 + a_2,$$

where we have set [see (3.16) and below (3.16) for the notation]

(3.24)
$$a_1 = P_x\left[\sum_{1 \leq k \leq k_*(A)} 1\{H_C < T_{\widetilde{B}'}\} \circ \theta_{R'_{k,A}} > \frac{K}{2}\right],$$

$$a_2 = P_x\left[\sum_{1 \leq k \leq k_*(G)} 1\{H_C < T_{\widetilde{B}'}; \pi_G(X_0) \notin A\} \circ \theta_{R'_k} > \frac{K}{2}\right],$$

and $R'_{1,A} = \inf\{R'_\ell; X_{R'_\ell} \in A \times \mathbb{Z}\}$, and for $k \geq 1$, $R'_{k+1,A} = \inf\{R'_\ell; R'_\ell > R'_{k,A}$ and $X_{R'_\ell} \in A \times \mathbb{Z}\}$, that is, $R'_{k,A}$, $k \geq 1$, stand for the successive times within $R'_\ell, \ell \geq 1$, when $\pi_G(X_{R'_\ell}) \in A$. There remains to bound $a_1$ and $a_2$. We first

write

$$a_1 \leq P_x\left[\sum_{2\leq k\leq k_*(A)} 1\{H_C < T_{\widetilde{B}'\cap(A\times\mathbb{Z})}\} \circ \theta_{R'_{k,A}} > \frac{1}{2}\left(\frac{K}{2}-1\right)\right]$$

(3.25)
$$+ P_x\left[\sum_{2\leq k\leq k_*(A)} 1\{T_{A\times\mathbb{Z}} < H_C < T_{\widetilde{B}'}\} \circ \theta_{R'_{k,A}} > \frac{1}{2}\left(\frac{K}{2}-1\right)\right]$$

$$\stackrel{\text{def}}{=} b_1 + b_2.$$

With the strong Markov property at times $R'_{k,A}$, recalling that when $k \geq 2$, $X_{R'_{k,A}} \in A \times \{-h', h'\}$, $P_x$-a.s., we find that for $\lambda > 0$,

(3.26)
$$b_1 \leq \exp\left\{-\frac{\lambda}{2}\left(\frac{K}{2}-1\right)\right\}$$
$$\times \left(\sup_{z\in A\times\{-h',h'\}} E_z[\exp\{\lambda 1_{\{H_C<T_{\widetilde{B}'\cap(A\times\mathbb{Z})}\}}\}]\right)^{k_*(A)}.$$

Choosing $\varepsilon = \frac{2\gamma-1}{2} \wedge \frac{\delta}{8}$ in Lemma 3.3, we see that for large $G$ in $\mathcal{G}_1$,

(3.27)
$$\sup_{z\in A\times\{-h',h'\}} P_z[H_C < T_{\widetilde{B}'\cap(A\times\mathbb{Z})}] \leq |C|h'^{-(2\gamma-1-\varepsilon)}.$$

Note also that for large $G$ in $\mathcal{G}_1$, with (3.22),

$$k_*(A)|C|h'^{-(2\gamma-1-\varepsilon)}$$

(3.28)
$$\stackrel{(3.16)}{\underset{(3.6)}{\leq}} c|C|\frac{|A|}{h'}|G|^{-15\delta/16}(\log|G|)^5 h'^{(1+\varepsilon)}|A|^{-1}|G|^{\delta/4}$$
$$\stackrel{(3.3)}{\leq} c|G|^{\delta/4+\delta/4-15\delta/16}h'^{\varepsilon}(\log|G|)^5$$
$$\stackrel{(3.6),(3.8)}{\leq} |G|^{-7\delta/16+\varepsilon} \leq |G|^{-\delta/4},$$

and as a result,

(3.29)
$$b_1 \leq \exp\left\{-\frac{\lambda}{2}\left(\frac{K}{2}-1\right) + k_*(A)|C|h'^{-(2\gamma-1-\varepsilon)}(e^\lambda-1)\right\}$$
$$\leq \exp\left\{-\frac{\lambda}{2}\left(\frac{K}{2}-1\right) + |G|^{-\delta/4}(e^\lambda-1)\right\}.$$

In an analogous fashion using in place of (3.27) the estimate

(3.30)
$$\sup_{z\in E} P_z[T_{A\times\mathbb{Z}} < H_C < T_{\widetilde{B}'}]$$
$$\stackrel{(3.5)}{\leq} h'|G|^{\delta/8-1} \stackrel{(3.6),(3.22)}{\leq} h'^{(1-2\gamma)} \leq |C|h'^{-(2\gamma-1-\varepsilon)},$$

we find in the case of $b_2$,

$$b_2 \leq \exp\left\{-\frac{\lambda}{2}\left(\frac{K}{2}-1\right) + |G|^{-\delta/4}(e^\lambda - 1)\right\}. \tag{3.31}$$

With similar arguments, we see that for large $G$ in $\mathcal{G}_1$, for $\lambda > 0$,

$$\begin{aligned}
a_2 &\leq \exp\left\{-\lambda\left(\frac{K}{2}-1\right) + k_*(G)h'|G|^{\delta/8-1}(e^\lambda-1)\right\} \\
&\leq \exp\left\{-\lambda\left(\frac{K}{2}-1\right) + |G|^{-\delta/4}(e^\lambda-1)\right\},
\end{aligned} \tag{3.32}$$

where we used $k_*(G)h'|G|^{\delta/8-1} \overset{(3.16)}{\leq} c|G|^{\delta/8-15\delta/16}(\log|G|)^5 \leq |G|^{-\delta/4}$.

Picking $\lambda = \frac{\delta}{4}\log|G|$ in (3.29), (3.31) and (3.32), we can choose a constant $K_1$ such that, for large $G$ in $\mathcal{G}_1$, when (3.22) holds,

$$\sup_{x \in E} P_x\left[\sum_{k \geq 1} 1\{\widetilde{R}_k < |G|^{2(1-\delta)}\} > K_1\right] \leq |G|^{-3}. \tag{3.33}$$

We then analyze the case where (3.22) is replaced with [cf. (3.6)]

$$h' = h. \tag{3.34}$$

For large $G$ in $\mathcal{G}_1$, under (3.34), for $K \geq 2, x \in E, \lambda > 0$, using analogous arguments as in (3.23) and the strong Markov property at time $R_m$, we obtain

$$\begin{aligned}
&P_x\left[\sum_{k \geq 1} 1\{\widetilde{R}_k < |G|^{2(1-\delta)}\} > K\right] \\
&\leq \exp\{-c(\log|G|)^5\} + P_x\left[\sum_{1 \leq m \leq m_*} 1\{H_C < T_{\widetilde{B}}\} \circ \theta_{R_m} > K\right] \\
&\leq \exp\{-\lambda(K-1)\}\left(\sup_{z \in G \times \{-h,h\}} E_z[\exp\{\lambda 1\{H_C < T_{\widetilde{B}}\}\}]\right)^{(m_*-1)}.
\end{aligned} \tag{3.35}$$

From (2.32) in Lemma 2.5, we know that, for large $G$,

$$\sup_{G \times \{-h,h\}} P_z[H_C < T_{\widetilde{B}}] \leq c|C|\frac{h}{|G|}. \tag{3.36}$$

Coming back to (3.35), we then find

$$\begin{aligned}
&P_x\left[\sum_{k \geq 1} 1\{\widetilde{R}_k < |G|^{2(1-\delta)}\} > K\right] \\
&\leq \exp\{-c(\log|G|)^5\} + \exp\left\{-\lambda(K-1) + cm_*|C|\frac{h}{|G|}(e^\lambda-1)\right\} \\
&\leq \exp\{-\lambda(K-1) + |G|^{-\delta/2}(e^\lambda-1)\},
\end{aligned}$$

where we used that, for large $G$ in $\mathcal{G}_1$, $m_* |C| \frac{h}{|G|} \stackrel{(3.13)}{\leq} 2|G|^{-\delta} (\log |G|)^5 |C| \stackrel{(3.3),(3.7)}{\leq} |G|^{-\delta/2}$. Picking $\lambda = \frac{\delta}{2} \log |G|$, we can choose a constant $K_2$ such that, for large $G$ in $\mathcal{G}_1$, when (3.34) holds,

$$(3.37) \quad \sup_{x \in E} P_x \left[ \sum_{k \geq 1} 1\{\widetilde{R}_k < |G|^{2(1-\delta)}\} > K_2 \right] \leq |G|^{-3}.$$

Combining (3.33) and (3.37), it now follows from translation invariance in the $\mathbb{Z}$-direction that with $K_0 = K_1 \vee K_2$, for large $G$ in $\mathcal{G}_1$,

$$\sup_{x \in E} P_x \left[ \sup_{u \in \mathbb{Z}} \sum_{k \geq 1} 1\{\widetilde{R}_k^u < |G|^{2(1-\delta)}\} > K_0 \right]$$

$$= \sup_{x \in G \times \{0\}} P_x \left[ \sup_{u \in \mathbb{Z}} \sum_{k \geq 1} 1\{\widetilde{R}_k^u < |G|^{2(1-\delta)}\} > K_0 \right]$$

$$\leq |G|^2 \sup_{x \in E} P_x \left[ \sum_{k \geq 1} 1\{\widetilde{R}_k < |G|^{2(1-\delta)}\} > K_0 \right] \leq |G|^{-1},$$

and this proves (3.12). $\square$

REMARK 3.4. The proof of Proposition 3.2 when (3.34) holds shows that for the class $\overline{\mathcal{G}}(d_0, \delta)$ of finite connected graphs satisfying (1.2) but not (3.1), if one defines, for $u$ in $\mathbb{Z}$ [see also (2.12)],

$$(3.38) \quad \overline{C}(u) = \overline{V} \times (u + [-w, w])$$

where $w \leq |\overline{V}| \wedge \left[\frac{h}{2}\right]$ and $|G|^{\delta/8} \geq |\overline{V}|$,

and introduces in analogy to (3.11), with $\widetilde{B}'(u)$ replaced by $\widetilde{B}(u)$, the successive returns to $\overline{C}(u)$ and departures from $\widetilde{B}(u)$, $\overline{R}_k^u, \overline{D}_k^u, k \geq 1$, one can find a positive $\overline{K}$ solely depending on $d_0$ and $\delta$ such that

$$(3.39) \quad \lim_{|G| \to \infty, G \in \overline{\mathcal{G}}} \sup_{x \in E} P_x \left[ \sup_{u \in \mathbb{Z}} \sum_{k \geq 1} 1\{\overline{R}^{u_k} < |G|^{2(1-\delta)}\} > \overline{K} \right] = 0.$$

This remark will be helpful in the next section when we derive a lower bound on the disconnection time for a large $G$ in $\mathcal{G}_{\text{ext}}$; see (3.10).

**4. Lower bound in presence of a transient pocket.** In this section we derive a lower bound on the disconnection time of a discrete cylinder when its base $G$ is large and contains a transient pocket [i.e., $\gamma > 1$, in the notation of (3.4)]. The basic result is Theorem 4.1; applications are given in Corollary 4.5, when $G$ is a truncated $r$-tree of depth $N$, or in Corollary 4.6 [see also (0.8)], when $G$ contains a ball of not too small radius modeled on an infinite graph where the heat kernel and volume growth conditions (0.5)(i) and (0.6) are fulfilled. Our methods also enable to derive a general lower bound on the disconnection time for large $G$ in $\mathcal{G}_{\text{ext}}$ (cf. Theorem 4.3), that is, when $\lambda_G$ is "close" to the extreme values in (1.9). Our key tools are Propositions 3.2 and 2.6. Our convention on constants in this section unless otherwise stated is the same as in Section 3; see above Remark 3.1. The definition of $\mathcal{G}_0$ appears at the beginning of Section 3 and corresponds to graphs $G$ where either $\lambda_G$ is "small" or a suitable "pocket" is present. Our main result in this section is the following:

THEOREM 4.1 (Transient pocket, $\gamma > 1$).

$$\lim_{|G|\to\infty, G\in\mathcal{G}_0} \sup_{x\in E} P_x[\mathcal{T}_G \leq |G|^{2(1-\delta)}] = 0. \tag{4.1}$$

PROOF. Choosing $\varepsilon_n = (\log n)^{-1}$ in (2.5), we see with (3.1) and (3.8) that

$$\lim_{|G|\to\infty, G\in\mathcal{G}_0\setminus\mathcal{G}_1} \sup_{x\in E} P_x[\mathcal{T}_G \leq |G|^{2(1-\delta)}] = 0. \tag{4.2}$$

We thus only need consider the case of a large $G$ in $\mathcal{G}_1$. We use the strategy outlined in Proposition 2.6. In the presence of a "transient pocket," that is, with $\gamma > 1$, we simply choose $M = M' = 1$, and $D_{1,1}(u) = D(u) = V \times (u + [0,v]) \subseteq C(u)$, for $u \in \mathbb{Z}$, in the notation of (2.41) and (3.7). So for a finite subset $S$ of $E$, "$S$ good for $D(u)$" is just the same as "$S$ thin in $D_{1,1}(u)(=D(u))$." The full strength of Proposition 2.6 will not be needed until Section 5. We also denote the image set of $X$, up to time $|G|^{2(1-\delta)}$ with

$$S = X_{[0,[|G|^{2(1-\delta)}]]}. \tag{4.3}$$

Using Propositions 2.6 and 3.2, our claim (4.1) will follow once we show that

$$\lim_{\substack{|G|\to\infty,\,x\in E \\ G\in\mathcal{G}_1}} \sup P_x\bigg[\text{for some } u\in\mathbb{Z}, S \text{ is not thin in } D(u), \text{ and for all} \tag{4.4}$$

$$u\in\mathbb{Z}, \sum_{k\geq 1} 1\{\widetilde{R}_k^u \leq |G|^{2(1-\delta)}\} \leq K_0\bigg]$$

$$= 0.$$

In order to contain the possible damage created by the few excursions reaching $C$, the next lemma will be useful.

LEMMA 4.2 ($\gamma > 1$).  *For large $G$ in $\mathcal{G}_1$ [cf. (3.2) and (3.7) for the notation],*

$$\sup_{x \in E} E_x \left[ \exp\left\{ \frac{c}{|V|^{1/\gamma}} |X_{[0, T_{\widetilde{B}'} - 1]} \cap C| \right\} \right] \leq 2. \tag{4.5}$$

PROOF.  Using a variation on Khashminskii's lemma [see also the proof of (2.33)], it suffices to show that

$$\sup_{x \in C} E_x [|X_{[0, T_{\widetilde{B}'} - 1]} \cap C|] \leq c|V|^{1/\gamma}. \tag{4.6}$$

To this end note that for large $G$ in $\mathcal{G}_1$, when $x = (g, w)$, one has

$$E_x[|X_{[0, T_{\widetilde{B}'} - 1]} \cap C|]$$

$$= \sum_{z \in C} P_x[H_z < T_{\widetilde{B}'}]$$

$$\leq |C| P_x[H_C \circ \theta_{T_{A \times \mathbb{Z}}} + T_{A \times \mathbb{Z}} < T_{\widetilde{B}'}]$$

$$+ E_x \left[ \sum_{n \geq 0} 1\{X_n \in C, n < T_{\widetilde{B}' \cap (A \times \mathbb{Z})}\} \right]$$

$$\stackrel{(3.5),(3.7)}{\leq} |C||h'||G|^{\delta/8 - 1} + c E_g^G \left[ \sum_{n \geq 0} 1\{Y_n \in V, n < T_A\} \right] \tag{4.7}$$

$$\stackrel{\substack{(3.3),(3.7) \\ (3.8),(3.10)}}{\leq} c|G|^{\delta/4} |G|^{1-\delta} (\log |G|)^3 |G|^{\delta/8 - 1} + c \sum_{k \geq 0} \left( \left( \frac{|V|}{k^\gamma} \right) \wedge 1 \right)$$

$$\leq 1 + c|V|^{1/\gamma} + c|V| \sum_{k > |V|^{1/\gamma}} k^{-\gamma}$$

$$\leq 1 + c|V|^{1/\gamma} \stackrel{(3.3)}{\leq} c|V|^{1/\gamma},$$

which proves (4.6). Our claim (4.5) follows.  □

We now prove (4.4). Considering $u_0 \in \mathbb{Z}$ with $|u_0| \leq |G|^2$, and $x \in E$, we find

$$P_x \left[ S \text{ is not thin in } D(u_0) \text{ and } \sum_{k \geq 1} 1\{\widetilde{R}_k^{u_0} \leq |G|^{2(1-\delta)}\} \leq K_0 \right]$$

$$\stackrel{(2.42)}{\leq} P_x\Big[ \text{ for some } k \leq K_0,$$

(4.8)
$$|\pi_G(X_{[0,T_{\widetilde{B}'(u_0)}-1]} \cap C(u_0))| \circ \theta_{\widetilde{R}_k^{u_0}} \geq \frac{|V|}{2K_0} \text{ or}$$

$$|\pi_{\mathbb{Z}}(X_{[0,T_{\widetilde{B}'(u_0)}-1]} \cap C(u_0))| \circ \theta_{\widetilde{R}_k^{u_0}} \geq \frac{v}{2K_0}\Big]$$

$$\text{with } v = |V| \wedge \left[\frac{h'}{2}\right]; \text{ see (3.7)}.$$

When $v = |V| < [\frac{h'}{2}]$, using Lemma 4.2 and the strong Markov at times $\widetilde{R}_k^{u_0}$, we find that the above probability is smaller than

(4.9)
$$P_x\Big[ \text{ for some } k \leq K_0, |X_{[0,T_{\widetilde{B}'(u_0)}-1]} \cap C(u_0)| \circ \theta_{\widetilde{R}_k^{u_0}} \geq \frac{|V|}{2K_0}\Big]$$

$$\leq 2K_0 e^{-c/(|V|^{1/\gamma})|V|/(2K_0)} \stackrel{(3.3)}{\leq} c\exp\{-c|G|^{\delta'(1-1/\gamma)}\}.$$

On the other hand, when $v = [\frac{h'}{2}] \leq |V|$, one either has [cf. (3.6)]

(4.10)  (i)  $|V| \geq \left[\frac{h'}{2}\right] = v$  and  $h' = [(\sqrt{|A|}|G|^{-\delta/8})^{1/\gamma}] < h,$

in which case, with a similar argument as above, the right-hand side of (4.8) is smaller than

$$2K_0\Big(\exp\Big\{-\frac{c}{2K_0}|V|^{1-1/\gamma}\Big\} + \exp\Big\{-c\Big(\frac{\sqrt{|A|}|G|^{-\delta/8}}{|V|}\Big)^{1/\gamma}\Big\}\Big)$$

$$\stackrel{(3.3),(1.9)}{\leq} c\exp\{-c|G|^{\delta'(1-1/\gamma)}\} + c\exp\{-c|G|^{7\delta/(32\gamma)}\},$$

otherwise, one has

(4.10)  (ii)  $|V| \geq \left[\frac{h'}{2}\right] = v$  and  $h' = h,$

in which case we instead use (2.33)(i) of Lemma 2.5 to bound the $\mathbb{Z}$-projection that appears in the last line of (4.8); we see that the right-hand side of (4.8) is smaller than

$$2K_0\Big(\exp\Big\{-\frac{c}{2K_0}|V|^{1-1/\gamma}\Big\} + \exp\Big\{-\frac{c}{\sqrt{t_G}}\frac{h}{5K_0}\Big\}\Big) \stackrel{(2.12)}{\leq} c\exp\{-c(\log|G|)^2\}.$$

Combining the above estimates, we see that for large $G$ in $\mathcal{G}_1$,

$$\sup_{x \in G \times \{0\}} P_x\Big[\text{for some } |u| \leq |G|^2, S \text{ is not thin in } D(u) \text{ and for all}$$

$$u \in \mathbb{Z}, \sum_{k \geq 1} 1\{\widetilde{R}_k^u \leq |G|^{2(1-\delta)}\} \leq K_0 \Bigg]$$

$$\leq c|G|^2 \exp\{-c(\log |G|)^2\}.$$

Using the fact that for $|u| > |G|^2$, $S$ is thin in $D(u)$ and translation invariance in the $\mathbb{Z}$-direction, we obtain (4.4). This concludes the proof of Theorem 4.1. □

The methods employed in the proof of Theorem 4.1 apply as well to the case of a large $G$ in $\mathcal{G}_{\text{ext}}$; see (3.10).

THEOREM 4.3.

(4.11) $$\lim_{|G| \to \infty, G \in \mathcal{G}_{\text{ext}}} \sup_{x \in E} P_x[\mathcal{T}_G \leq |G|^{2(1-\delta)}] = 0.$$

PROOF. Using (4.2), we see that we can replace $\mathcal{G}_{\text{ext}}$ in (4.11) with

$$\widetilde{\mathcal{G}}_{\text{ext}} = \{G \in \mathcal{G}_{\text{ext}}; \lambda_G \geq |G|^{-\delta/10}\}.$$

We now employ an analogous strategy as explained below (4.2). We choose again $M = M' = 1$, $D_{1,1}(u) = D(u)$, where, for $u \in \mathbb{Z}$,

(4.12) $$D(u) = V \times \left(u + \left[0, |V| \wedge \left[\frac{h}{2}\right]\right]\right),$$

where $V$ is a connected subset of $G$ with $|V| = [|G|^{\delta/8}]$, and $h$ as in (2.12). Note that for large $G$ in $\widetilde{\mathcal{G}}_{\text{ext}}$, $\frac{h}{2} \leq |V|$. Using (3.39), we only need to show, with $S$ as in (4.3),

(4.13) $$\lim_{|G| \to \infty, G \in \widetilde{\mathcal{G}}_{\text{ext}}} \sup_{x \in E} P_x\Bigg[\text{for some } u \in \mathbb{Z}, S \text{ is not thin in } D(u) \text{ and}$$
$$\text{for all } u \in \mathbb{Z}, \sum_{k \geq 1} 1\{\overline{R}_k^u \leq |G|^{2(1-\delta)}\} \leq \overline{K}\Bigg] = 0.$$

Note that for large $G$ in $\widetilde{\mathcal{G}}_{\text{ext}}$, $|u_0| \leq |G|^2$, and $x \in E$, one has

$$P_x\Bigg[S \text{ is not thin in } D(u_0) \text{ and } \sum_{k \geq 1} 1\{\overline{R}_k^{u_0} \leq |G|^{2(1-\delta)}\} \leq \overline{K}\Bigg]$$

$$\stackrel{(2.33)}{\leq} 2\overline{K}\left(\exp\left\{-\frac{c}{t_G}\frac{|V|}{2\overline{K}}\right\} + \exp\left\{\frac{c}{\sqrt{t_G}}\frac{1}{2\overline{K}}\left[\frac{h}{2}\right]\right\}\right)$$

$$\leq c\exp\{-c(\log|G|)^2\},$$

where the constants matter-of-factly only depend on $d_0$ and $\delta$ (as in Section 2). We can then conclude the proof of Theorem 4.3 as we did below (4.10)(ii). □

REMARK 4.4. (i) $\mathcal{G}_{\text{ext}}(d_0, \delta)$ contains, on the one hand, large, "one-dimensional" finite graphs such as $\mathbb{Z}/N\mathbb{Z}$ or $[0, N]$ ($d_0 \geq 2$ and $0 < \delta < \frac{1}{2}$, arbitrary), for which $\lambda_G$ is of order $|G|^{-2}$ and on the other hand, for $d_0 \geq 3$, $0 < \delta < \frac{1}{2}$, large expanders (cf. [13]) for which $\lambda_G$ is order a positive constant.

(ii) As an immediate consequence of Theorems 1.2 and 4.3, we thus see that, for $d_0 \geq 2$, $0 < \delta < \frac{1}{2}$, $\varepsilon > 0$,

$$(4.14) \quad \lim_{|G| \to \infty, G \in \mathcal{G}_{\text{ext}}(d_0, \delta)} \inf_{x \in E} P_x[|G|^{2(1-\delta)} \leq \mathcal{T}_G \leq |G|^2 (\log |G|)^{4+\varepsilon}] = 1.$$

This, of course, immediately implies that (0.9) holds for sequences satisfying (1.2) and (0.10).

We will now describe two applications of Theorem 4.1. The first application concerns the case where, for $N \geq 1$,

(4.15)    $G_N$ is a rooted $r$-tree of depth $N$, with root denoted by $g_*$;

here $r \geq 2$ is an integer (and the case $r = 2$ corresponds to the rooted binary tree of depth $N$). Clearly, one has $|G_N| = 1 + r \ldots r^N = r^{N+1} - 1$. We set $d_0 = r + 1$, and write $E_N = G_N \times \mathbb{Z}$.

COROLLARY 4.5 [*Under* (4.15)].

$$(4.16) \quad \lim_{N \to \infty} \inf_{x \in E_N} P_x[|G_N|^{2(1-\delta)} \leq \mathcal{T}_{G_N} \leq |G_N|^2 (\log |G_N|)^{4+\varepsilon}] = 1$$

*for* $\varepsilon, \delta > 0$.

PROOF. The upper bound follows from Theorem 1.2. For the lower bound without loss of generality, we choose $\delta \in (0, \frac{1}{2})$. We view $G_N$ as a subset, namely, the open ball with center $g_*$ and radius $N + 1$, of $G_\infty$ the infinite $r$-tree with root $g_*$. In the notation of (3.2) we pick $A_N$ and $V_N$ as open balls with center $g_*$:

(4.17)    $A_N = B(g_*, N) \supseteq V_N = B(g_*, \rho_N)$

(4.18)    where $\rho_N$ is an integer such that $|G_N|^{\delta/16} \geq |V_N| \geq \dfrac{1}{r+1} |G_N|^{\delta/16}$.

Note that the random walks on $G_N$ or on $G_\infty$ killed when exiting $A_N$ do agree. Moreover, with hopefully obvious notation (see also the beginning of Section 1), $d(Y_k, g_*)$ under $P_g^{G_\infty}$ is distributed as a simple random walk on the nonnegative integers reflected at 0, with jump probability "to the right" equal to $\frac{r}{r+1} > \frac{1}{2}$. It thus follows that, for $N \geq 1$, $g \in A_N$, $k \geq 0$,

(4.19) $\quad P_g^{G_N}[X_k = g, k < T_{A_N}] \leq P_g^{G_\infty}[X_k = g] \leq e^{-\mu(r)k} \qquad$ with $\mu(r) > 0$,

using a comparison with the simple random walk on $\mathbb{Z}$, jumping to the right with probability $\frac{r}{r+1}$, for the last inequality. In view of Theorem 4.1, the claim (4.16) will follow once we show that

(4.20) $\quad$ for large $N$, $\quad G_N \in \mathcal{G}_0\bigg(d_0 = r+1, \delta, \delta' = \frac{\delta}{16}, \gamma = \frac{100}{\delta},$

$$a = \sup_{k \geq 1} k^{100/\delta} e^{-\mu(r)k}\bigg).$$

Since $|A_N|/|G_N|$ remains bounded away from 0, it follows that for large $N$ either (3.1) or (3.3) holds; see also (1.9). Although we do not explicitly need the following fact, it is of interest to remark that with (59) and (60) in Chapter 5 of [2], $\lambda_{G_N}$ is of order $|G_N|^{-1}$ for large $N$ ([2] discusses the spectral gap attached to the discrete time walk, which can be compared to $\lambda_{G_N}$ by a bounded multiplicative factor depending on $d_0 = r+1$). Since $\delta < \frac{1}{2}$, in fact, (3.1) does not hold for large $N$. Clearly, in view of (4.19) and the choice of $a$ in (4.20), (3.4) holds as well. As a result, (4.20) and, hence, our claim (4.16) will follow once we show that

(4.21) $\qquad\qquad\qquad$ for large $N$, (3.5) holds.

To this end, note first that, for $g \in G_N \setminus A_N$, that is, when $d(g, g_*) = N$,

$$P_g^{G_N}[H_{V_N} \circ \theta_1 < T_{A_N} \circ \theta_1] \leq r^{-(N-\rho_N)}.$$

Indeed, the walk on $G_N$ and $G_\infty$ coincide up to the exit time from $A_N$, and the distance to $g_*$ of the walk on $G_\infty$ has the law described below (4.19). The above bound now follows from the application of the simple Markov property and standard estimates for the biased simple random walk on $\mathbb{Z}$ [note that the ball defining $V_N$ in (4.17) is open]. It thus follows that, for $T > 0$, $N \geq 1$,

(4.22) $\qquad\qquad \sup_{g \in A_N^c} P_g^{G_N}[H_{V_N} < T] \leq T r^{-(N-\rho_N)}.$

Then observe that in the notation of (3.6), with (4.17),

$$\sup_{N \geq 1} h'_N r^{-(N/\gamma)(1/2 - \delta/8)} \leq \nu(r, \delta) < \infty,$$

and hence, there is for large $N$, for any $x \in \widetilde{B}'_N$, a probability at least $p(r,\delta) > 0$ to exit $\widetilde{B}'_N$ before time $r^{(N/\gamma)(1-\delta/4)}$ under $P_x$. It now follows that, for large $N$, and $x \in A_N^c \times \mathbb{Z}$,

$$P_x[H_{C_N} < T_{\widetilde{B}'_N}] \leq P_x[H_{C_N} < T] + P_x[T < T_{\widetilde{B}'_n}]$$
(4.23)
$$\leq Tr^{-(N-\rho_N)} + (1-p)^{[T/r^{N/\gamma(1-\delta/4)}]},$$

for $T$ a positive integer, using (4.22) for the first term in the last line and the remark above (4.23). Choosing $T = [r^{N/\gamma(1-\delta/4)+N\delta/100}]$, and noting in view of (4.18) that $\rho_N \sim N\delta/16$, and $1/\gamma(1-\delta/4) + \delta/100 + \delta/16 < \delta/8$ (recall $\gamma = 100/\delta$), it follows that for large $N$ the expression in the second line of (4.23) is smaller than $|G_N|^{\delta/8-1}$. This is more than enough to prove (4.21). This concludes the proof of (4.16). □

We now turn to the second application of Theorem 4.1. We consider an infinite connected graph $G_\infty$ with degree bounded by $d_0 \geq 2$, with polynomial volume growth (0.6) and such that the random walk satisfies the upper bound (0.5)(i) [we do not require (0.5)(ii)]. We assume in this section that

(4.24) $$\alpha > \beta \geq 2,$$

and in view of (0.5)(i), the walk is transient on $G_\infty$. When $\beta > 2$, the walk is sub-diffusive [for instance, with (0.5)(i) and (0.6), the expected distance from the starting point at time $k$ is uniformly bounded by const $k^{1/\beta}$], a feature often referred to as anomalous diffusion. For a thorough investigation of such walks and examples, we refer to [3], [10] and [11].

We now consider a sequence of connected finite graphs $G_N, N \geq 1$, with degree bounded by $d_0 \geq 2$, such that $\lim_N |G_N| = \infty$, and for large $N$, there is $r_N > 1$, and $g_N \in G_N$, such that

$$B(g_N, r_N) \subseteq G_N \text{ is isomorphic to some open ball of radius } r_N \text{ in } G_\infty$$
(4.25)

[i.e., there is a bijection between $B(g_N, r_N)$ and an open ball of radius $r_N$ in $G_\infty$, which preserves the degree, and such that pairs of points in $B(g_N, r_N)$ are neighbors if and only if their images in $G_\infty$ are neighbors], and for a suitable $\eta \in (0,1)$, and sequence $\varphi_n$ such that $\varphi_n = o(n^\varepsilon)$, for each $\varepsilon > 0$,

(4.26) $$|B(g_N, r_N)| \geq \min(|G_N|\varphi_{|G_N|}^{-1}, \lambda_{G_N}^{-1/2}|G_N|^\eta) \qquad \text{for large } N.$$

COROLLARY 4.6. *Under the above assumptions, for all $\delta > 0, \varepsilon > 0$,*

(4.27) $$\lim_{N\to\infty} \inf_{x\in E_N} P_x[|G_N|^{2(1-\delta)} \leq \mathcal{T}_{G_N} \leq |G_N|^2(\log|G_N|)^{4+\varepsilon}] = 1.$$

PROOF. As for Corollary 4.5, only the lower bound on $\mathcal{T}_{G_N}$ is of concern. We choose $0 < \delta < \frac{1}{2} \wedge \eta$, see (4.26). For the remainder of the proof, all constants $c$ may depend on $d_0, \delta, \alpha, \beta, \kappa_i, \widetilde{\kappa}_i$, $i = 1, 2$ [cf. (0.5), (0.6)], $\eta$.

We will apply Theorem 4.1, and to this end, choose for large $N$

$$A_N = B(g_N, r_N), \quad V_N = B(g_N, \rho_N) \quad (4.28)$$

$$\text{with } |G_N|^{\delta/8} \geq |V_N| \geq \frac{1}{d_0}|G_N|^{\delta/8}.$$

The claim (4.26) will follow once we show that

(4.29) $\quad$ for large $N$, $\quad G_N \in \mathcal{G}_0\left(d_0, \delta, \delta' = \frac{\delta}{16}, \gamma = \frac{\alpha}{\beta}, a = \kappa_1\right).$

Given that $\delta < \eta$, in view of (4.26), (4.25) and (0.5)(i), we only need to prove that

(4.30) $\quad$ for large $N$, (3.5) holds.

With (4.26), (0.5)(i) and (0.6), we observe that for $g \in G_N$ with $d(g, g_N) = [\frac{r_N}{4}]$ and $T \geq 1$,

$$P_g^{G_N}\left[\sup_{1 \leq k \leq T} d(g, Y_k) \geq \frac{r_N}{4}\right] \leq \sum_{1 \leq k \leq T} \frac{c}{k^{\alpha/\beta}} r_N^\alpha \exp\left\{-c\left(\frac{r_N^\beta}{k}\right)^{1/(\beta-1)}\right\}$$

(4.31)

$$\leq c r_N^\alpha \exp\left\{-c\left(\frac{r_N^\beta}{T}\right)^{1/(\beta-1)}\right\},$$

using (4.24) in the last step. With (4.26), (4.28) and (0.6), we see that, for large $N$ (recall $\eta > \delta$),

(4.32) $\quad \frac{r_N}{8} \geq c|G_N|^{\eta/\alpha} \geq c'|G_N|^{\delta/(8\alpha)} \geq \rho_N,$

and with (3.6),

$$h_N' \leq (|B(g_N, r_N)|^{1/2}|G_N|^{-\delta/8})^{1/\gamma} \leq c r_N^{\beta/2}|G_N|^{-\delta/(8\gamma)}$$

(4.33)

$$\text{where } \lim_N h_N' = \infty.$$

The same argument as for (4.23) shows that, for large $N$, and $x \in A_N^c \times \mathbb{Z}$,

$P_x[H_{C_N} < T_{\widetilde{B}_N'}]$

$$\leq P_x[H_{C_N} < T] + P_x[T < T_{\widetilde{B}_N'}]$$

(4.34)

$$\overset{(4.31)}{\leq} c r_N^\alpha \exp\left\{-c\left(\frac{r_N^\beta}{T}\right)^{1/(\beta-1)}\right\} + (1-c)^{[T/(c r_N^\beta |G_N|^{-\delta/(4\gamma)})]}$$

for $T$ a positive integer.

Choosing $T = [r_N^\beta |G_N|^{-\delta/(8\gamma)}]$ and observing that $r_N^\alpha \leq c|G_N|$, it follows that for large $N$ the last line of (4.34) is smaller than $|G_N|^{-1}$. This is more than enough to show (4.30), and concludes the proof of Corollary 4.6. □

REMARK 4.7. (i) Corollary 4.6 applies, in particular, when $G_N$ is some ball of radius $N$ in $G_\infty$ (with possibly variable center), with $G_\infty$ as specified above Corollary 4.6.

(ii) One may also obtain (4.27) in situations where Corollary 4.6 does not directly apply. For instance, assume that $G_N$ satisfies the assumptions of Corollary 4.6 and $\widetilde{G}_N$ is a sequence of connected graphs with degree bounded by $\widetilde{d}_0$, such that, for some $\widetilde{\eta} > 0$, and large $N$,

$$|\widetilde{G}_N| \leq |G_N|^{\widetilde{\eta}} \tag{4.35}$$

and

$$\lambda_{\widetilde{G}_N} \geq \lambda_{G_N}, \tag{4.36}$$

then

(4.37)    (4.27) holds true    with $\overline{G}_N = G_N \times \widetilde{G}_N$ in place of $G_N$.

Indeed, we simply choose $\overline{A}_N = A_N \times \widetilde{G}_N$, $\overline{V}_N = V_N \times \{\widetilde{x}_N\}$, with $\widetilde{x}_N$ some point in $\widetilde{G}_N$. Using the fact that $\lambda_{\overline{G}_N} = \min(\lambda_{G_N}, \lambda_{\widetilde{G}_N})$, a product formula in the spirit of (1.5) to gain control over the random walk on $\overline{G}_N$, and (4.34), one sees that for $\delta$ as in the proof of Corollary 4.6 with a suitably large enough $a$,

$$\text{for large } N, \qquad G_N \in \mathcal{G}_0\left(d_0 + \widetilde{d}_0, \delta, \delta' = \frac{\delta}{16(1+\widetilde{\eta})}, \gamma = \frac{\alpha}{\beta}, a\right). \tag{4.38}$$

(iii) We can choose $G_\infty$ to be $\mathbb{Z}^d$ $d, \geq 3$, with its usual graph structure, so that $\alpha = d > \beta = 2$; see Remark 1.2 of [10]. We can then apply Corollary 4.6, when $G_N = (\mathbb{Z}/N\mathbb{Z})^d$, $N \geq 1$, and recover Theorems 1.1 and 2.1 of [8], when $d \geq 3$. The case $d = 1$ is covered by Remark 4.4(2). The case $d = 2$ will follow from the results in the next section. Periodic boundary conditions play no role here, and Corollary 4.6 applies just as well to the case $G_N = [0, N]^d$, $N \geq 1$; see (4.25) and (4.26).

(iv) When $G_N = [0, [N^\lambda]]^{d-1} \times [0, N]$, $N \geq 1$, with $\lambda \in (0, \infty)$, $d \geq 3$, then for large $N$, $\lambda_{G_N} N^{2(\lambda \vee 1)}$ remains bounded and bounded away from zero. Choosing $A_N$ to be a ball in $G_N$ with radius small multiple of $N^{\lambda \wedge 1}$, and suitable center, Corollary 4.6 applies [cf. (4.26)] if $d(\lambda \wedge 1) > \lambda \vee 1$, that is, when $\lambda \in (\frac{1}{d}, d)$. Further, when $d \geq 4$, $\lambda > 1$, Remark 4.7(2) applies (with $\widetilde{G}_N = [0, N]$), and we see that (4.27) holds true when $\lambda > \frac{1}{d}$.

**5. Lower bound in the presence of recurrent pockets.** We now derive a lower bound on the disconnection time of a discrete cylinder that applies to cases where $G$ is large and contains a recurrent pocket. This is substantially more delicate to handle than the case of transient pockets treated in Section 4. In particular, we make full use of Proposition 2.6 when the notion "$S$ is good in $D(u)$" [cf. (2.43)] involves a kind of renormalization step. We have shown in Proposition 3.2 the rarefaction of excursions between $C$ and the complement of $\widetilde{B}'$, taking place before time $|G|^{2(1-\delta)}$, when $G$ is large in $\mathcal{G}_0$ for arbitrary $\gamma > \frac{1}{2}$. However, we are unable to extend Theorem 4.1 to the case $\frac{1}{2} < \gamma \leq 1$. We need additional assumptions to tame the possible recurrence properties of the walk on $G$. The main result is Theorem 5.2; applications are then discussed in Corollary 5.3 and Remark 5.5.

Assuming $\frac{1}{2} < \gamma \leq 1$, we now describe the sub-class of $\mathcal{G}_0$ (cf. beginning of Section 3) consisting of $G$ in $\mathcal{G}_0$ such that, when (3.1) does not hold,

(5.1) $V$ in (3.1) is a geodesic segment [i.e., $V = \{g_\ell; 0 \leq \ell < |V|\}$, with $d(g_\ell, g_{\ell'}) = |\ell - \ell'|$, for $0 \leq \ell, \ell' < |V|$],

for $W \subseteq V$, a geodesic segment of length $m \geq 2$, and $J \subseteq \mathbb{Z}$ an interval of length $[m^{\beta/2}(\log m)^{-\beta/2}]$,

(5.2) $P_x[H_{W \times J} < T_{A \times \mathbb{Z}}] \leq a' \max\left(1, \frac{d(g,W)}{m}, \frac{d_{\mathbb{Z}}(u,J)^{2/\beta}}{m}\right)^{-\nu},$

for $x = (g,u) \in E$, with $d(g,W) = \inf\{d(g,g'); g' \in W\}$, and $d_{\mathbb{Z}}(u,J)$ analogously defined,

(5.3) for $W, J, m$, as above, and $\widetilde{g} \in W$,
$E_x[|\pi_G(X_{[0,T_{(A \cap B(\widetilde{g}, m \log m)) \times \mathbb{Z}}]}) \cap (W \times J))|] \leq a_G m (\log m)^{-\mu},$
for $x \in E$ and when $\pi_{\mathbb{Z}}$ replaces $\pi_G$, the right-hand side is replaced with $a_{\mathbb{Z}} m^{\beta/2} (\log m)^{-\mu - \beta/2}$.

Our assumptions on the parameters that appear above are

(5.4) $\quad \beta \geq 2, \quad \mu > 0, \quad \mu + \nu > 1, \quad a', \ a_G, a_{\mathbb{Z}} \geq 1.$

We denote with $\widetilde{\mathcal{G}}_0(d_0, \delta, \delta', \gamma, \mu, \nu, \beta, a, a', a_G, a_{\mathbb{Z}})$ the above defined class (we recall that here $\frac{1}{2} < \gamma \leq 1$). Unless otherwise stated, for the remainder of this section $c$ denotes a positive constant possibly depending on the above parameters, with a corresponding meaning for the expression "for large $G$ in ...."

REMARK 5.1. Let us give a word of comment on the above class $\widetilde{\mathcal{G}}_0$. The parameter $\beta \geq 2$ has a similar interpretation as in (0.5), with $\beta > 2$, enabling "anomalous diffusion" in the pocket $A$ (and hence, much faster displacements of the $\mathbb{Z}$-component than the $G$-component). The most restrictive assumption is (5.3), with $\mu > 0$. It rules out applications to pockets modeled on a suitable ball in an infinite graph of bounded degree satisfying (0.5), when $1 + \frac{\beta}{2} > \alpha$, which are instances of so-called "very strong recurrence" since $\beta > \alpha$; see Proposition 3 of Barlow [3]. In Corollary 5.3 we consider the recurrent situation $\beta \geq \alpha \geq (1 + \frac{\beta}{2}) \vee (\beta - 1)$, with $\beta \geq 2$, and can choose $\nu = \alpha - \frac{\beta}{2}$ in (5.2), and $\mu$ arbitrarily close to 1 in (5.3). The most delicate situation arises when $\alpha = 1 + \frac{\beta}{2}$.

THEOREM 5.2.

$$\lim_{|G| \to \infty, G \in \widetilde{\mathcal{G}}_0} \sup_{x \in E} P_x[\mathcal{T}_G \leq |G|^{2(1-\delta)}] = 0. \tag{5.5}$$

PROOF. We begin with some preparatory remarks. With Theorem 4.3 [see also (3.10)] our claim follows once we prove (5.5) with $\widetilde{\mathcal{G}}_0$ replaced by

$$\widetilde{\mathcal{G}} = \{G \in \widetilde{\mathcal{G}}_0; |G|^{-\delta/10} > \lambda_G > |G|^{-2(1-\delta)}(\log |G|)^{-1}\}. \tag{5.6}$$

Further, reducing $\nu$ and $V$ if necessary [cf. (5.2) and (3.3)], we can assume that

$$\delta' < \frac{\delta}{10}, \qquad |G|^{\delta'} \geq |V| \geq \frac{1}{d_0}|G|^{\delta'} \quad \text{and} \quad 0 < \nu < 1 \tag{5.7}$$

with $\mu + \nu > 1$.

We then choose $\rho(\mu, \nu) \in (0, 1)$ so that

$$\mu + \rho\nu > 1. \tag{5.8}$$

For large $G$ in $\widetilde{\mathcal{G}}$, we introduce the integers $M, M', L$ such that

$$M = [(\log|G|)^\rho], \qquad \tfrac{1}{2}|V| \leq M^{\beta/2} L^{\beta/2} < |V| \quad \text{and} \tag{5.9}$$
$$M' = [M^{\beta/2}(\log L)^{\beta/2}].$$

Recall that $V$ is a geodesic segment [cf. (5.1)] and in the notation of (2.39)–(2.41), we define

$$V_i = \{g_{(i-1)L}, \ldots, g_{iL}\}, \qquad J_j = [(j-1)L', jL'] \tag{5.10}$$

with $L' = [L^{\beta/2}(\log L)^{-\beta/2}]$.

Let us mention that when $\beta > 2$, the sets $D_{i,j} = V_i \times J_j$ are thin vertically elongated "rectangles" and their union $D$ [see notation below (2.41)] is also a thin vertically elongated rectangle contained in the rectangle $C$, with height $M'L'$ comparable to the height of $C$; see (3.7) and (5.12) below. Using similar arguments as for Theorem 4.1 (see, in particular, the end of the proof) our claim will follow once we show that

$$\text{(5.11)} \quad \lim_{|G|\to\infty, G\in\widetilde{\mathcal{G}}} |G|^2 \sup_{x\in E} P_x\left[S \text{ is not good in } D \text{ and } \sum_{k\geq 1} 1\{\widetilde{R}_k \leq |G|^{2(1-\delta)}\} \leq K_0\right] = 0,$$

with the notation of (4.3), (2.41), (2.43) and (3.11).

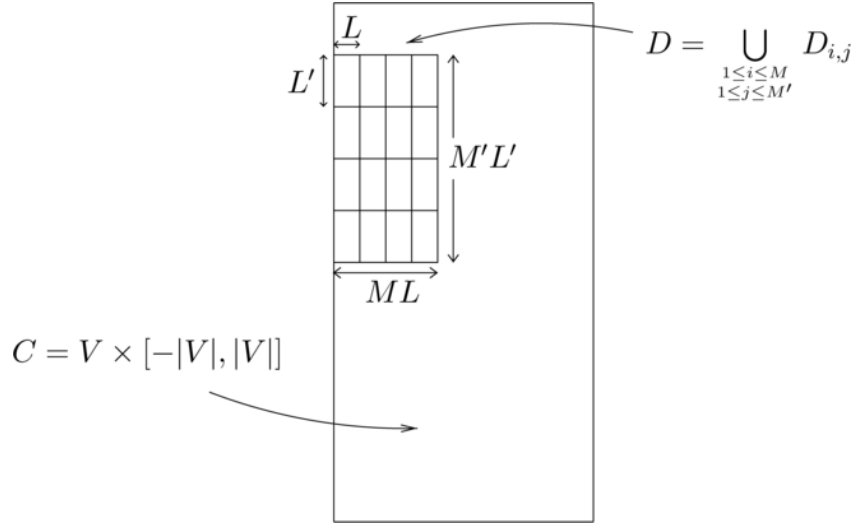

FIG. 3. *An illustration of the sets $D_{i,j}$, $D$ and $C$, when (5.12) holds.*

As a last reduction note that in view of (5.6), (5.7), (3.3), (3.6), here ($\frac{1}{2} < \gamma \leq 1$), for large $G$ in $\widetilde{\mathcal{G}}$, one has

$$\text{(5.12)} \quad |V| < \frac{h'}{2}.$$

We will now bound the probability that appears in (5.11). On the event inside the probability in this display, we can find a subset of $\{1,\ldots,M\} \times \{1,\ldots,M'\}$, consisting of ordered pairs $(i,j)$ such that $S$ is not thin in $D_{i,j}$, and either this subset consists of $[\frac{M}{2}]$ elements having different first coordinates, or this subset consists of $M'$ elements having different second coordinates; see (2.43). In the first case, we denote with $\mathcal{H}_0$ this subset and observe that for each $(i,j)$ in $\mathcal{H}_0$, one has [with the notation below (2.12)]

$$\text{(5.13)} \quad \sum_{1\leq k\leq K_0} |\pi_G(X_{[0,T_{\widetilde{B}'}]} \cap D_{i,j})| \circ \theta_{\widetilde{R}_k} \geq L/2, \quad \text{or}$$

$$\sum_{1 \leq k \leq K_0} |\pi_{\mathbb{Z}}(X_{[0,T_{\widetilde{B}'}]} \cap D_{i,j})| \circ \theta_{\widetilde{R}_k} \geq [L^{\beta/2}(\log L)^{-\beta/2}]/2.$$

In the second case, we denote with $\mathcal{H}_1$ the set of the $M$ elements in the above subset having second coordinate of the form $[\ell^{\beta/2}(\log L)^{\beta/2}]$, for $1 \leq \ell \leq M$, so that now (5.13) holds for all $(i,j) \in \mathcal{H}_1$.

With a rough counting argument, there are at most $2^M (M')^{M/2}$ possible choices for $\mathcal{H}_0$ and $M^M$ possible choices for $\mathcal{H}_1$. Thus, for large $G$ in $\widetilde{\mathcal{G}}$, with (5.9),

(5.14)    there are at most $\exp\{c(\log |G|)^\rho \log \log |G|\}$ possible choices for $\mathcal{H}_0$ or $\mathcal{H}_1$.

We will now bound the probability that, for each $(i,j)$ in $\mathcal{H}_0$ or $\mathcal{H}_1$, (5.13) happens. To this end, we consider for large $G$ in $\widetilde{\mathcal{G}}$ some $\mathcal{H}_0$ as above, and the "vertical" or "horizontal" segments of the form

(5.15)    $$U = \{g\} \times J_k \subseteq \bigcup_{\mathcal{H}_0} D_{i,j}, \qquad W = V_\ell \times \{u\} \subseteq \bigcup_{\mathcal{H}_0} D_{i,j}.$$

One then has for $z \in E$, with hopefully obvious notation,

$$E_z\left[\sum_{\mathcal{H}_0} |\pi_G(X_{[0,T_{\widetilde{B}'}]} \cap D_{i,j})|\right]$$

$$= E_z\left[\sum_U 1\{H_U < T_{\widetilde{B}'}\}\right]$$

(5.16)   $$\leq \sum_U P_z[T_{A \times \mathbb{Z}} < H_U < T_{\widetilde{B}'}] + \sum_{(i,j) \in \mathcal{H}_0} \sum_{U \subseteq D_{i,j}} P_z[H_U < T_{\widetilde{B}' \cap (A \times \mathbb{Z})}]$$

$$\leq |C| P_z[T_{A \times \mathbb{Z}} < H_C < T_{\widetilde{B}'}]$$

$$+ \sum_{(i,j) \in \mathcal{H}_0} E_z\bigg[H_{D_{i,j}} < T_{\widetilde{B}' \cap (A \times \mathbb{Z})},$$

$$\sum_{U \subseteq D_{i,j}} P_{X_{H_{D_{i,j}}}}[H_U < T_{\widetilde{B}' \cap (A \times \mathbb{Z})}]\bigg]$$

$$\overset{(3.5)}{\leq} |C| h' |G|^{\delta/8 - 1} + Q \overset{(5.6),(5.7)}{\leq} 1 + Q,$$

where $Q$ stands for the expression in the fourth line of (5.16). We now bound $Q$. First note that for $(i,j) \in \mathcal{H}_0$ and for $x \in D_{i,j}$ (playing the role of $X_{H_{D_{i,j}}}$ in the expectation entering $Q$), one has

$$
\begin{aligned}
(5.17) \quad & \sum_{U \subseteq D_{i,j}} P_x[H_U < T_{\widetilde{B}' \cap (A \times \mathbb{Z})}] \\
& \leq E_x\left[\sum_{k \geq 1} |\pi_G(X_{[0,T_{(B(g_{iL}, L \log L) \cap A) \times \mathbb{Z}]}} \cap D_{i,j})| \circ \rho_k, \rho_k < T_{\widetilde{B}' \cap (A \times \mathbb{Z})}\right],
\end{aligned}
$$

where analogously to (1.12), $\rho_k, k \geq 1$, are the successive return times to $D_{i,j}$ after leaving $(B(g_{iL}, L \log L) \cap A) \times \mathbb{Z}$. Using (5.2) with $m = L$, we find that

$$
(5.18) \quad P_{x'}[H_{D_{i,j}} < T_{A \times \mathbb{Z}}] \leq c(\log L)^{-\nu} \leq \tfrac{1}{2}
$$

for $x' \in (B(g_{iL}, L \log L)^c \cap A) \times \mathbb{Z}$.

Then using the strong Markov property at times $\rho_k$, we see that

$$
\begin{aligned}
\sum_{U \subseteq D_{i,j}} & P_x[H_U < T_{\widetilde{B}' \cap (A \times \mathbb{Z})}] \\
& \leq \left(\sum_{k \geq 1} \frac{1}{2^{k-1}}\right) \sup_{\widetilde{x} \in D_{i,j}} E_{\widetilde{x}}[|\pi_G(X_{[0, T_{(B(g_{iL}, L \log L) \cap A) \times \mathbb{Z}]}} \cap D_{i,j})|] \\
& \stackrel{(5.3)}{\leq} 2a_G L(\log L)^{-\mu}.
\end{aligned}
$$

Coming back to the fourth line of (5.16), we obtain for large $G$ in $\widetilde{\mathcal{G}}$, $z \in E$, $\mathcal{H}_0$ as above (5.13),

$$
\begin{aligned}
Q & \leq cL(\log L)^{-\mu} \sum_{(i,j) \in \mathcal{H}_0} P_z[H_{D_{i,j}} < T_{\widetilde{B}' \cap (A \times \mathbb{Z})}] \\
(5.19) \quad & \leq cL(\log L)^{-\mu} \sup_{\widetilde{z} \in D} \sum_{(i,j) \in \mathcal{H}_0} P_{\widetilde{z}}[H_{D_{i,j}} < T_{\widetilde{B}' \cap (A \times \mathbb{Z})}] \\
& \stackrel{(5.2)}{\leq} cL(\log L)^{-\mu}\left(1 + \sum_{1 \leq k \leq M} \frac{c}{k^\nu}\right) \leq cL(\log L)^{-\mu} M^{1-\nu},
\end{aligned}
$$

where we used the strong Markov property at time $H_D$ [cf. (2.41)] for the first inequality, and the structure of $\mathcal{H}_0$ [see above (5.13)] together with the fact that $V$ is a geodesic segment, for the second inequality. Coming back to (5.16), we see that for large $G$ in $\widetilde{\mathcal{G}}$ and any $\mathcal{H}_0$ as above (5.13),

$$
(5.20) \quad \sup_{z \in E} E_z\left[\sum_{(i,j) \in \mathcal{H}_0} |\pi_G(X_{[0, T_{\widetilde{B}'}]} \cap D_{i,j})|\right] \leq cL(\log L)^{-\mu} M^{1-\nu}.
$$

A similar bound holds as well with $\mathcal{H}_1$ in place of $\mathcal{H}_0$. Indeed, one simply needs to replace the sum in the last line of (5.19) with $\sum_{\ell=1}^M c(\ell^{\beta/2})^{-(2\nu)/\beta} \leq cM^{1-\nu}$; see (5.2) and below (5.13).

We now want to derive similar controls to (5.20) when $\pi_{\mathbb{Z}}$ replaces $\pi_G$. With analogous arguments as for (5.17), we see that, for $x \in D_{i,j}$ and $W$ as in (5.15),

$$\sum_{W \subseteq D_{i,j}} P_x[H_W < T_{\widetilde{B}' \cap (A \times \mathbb{Z})}]$$

(5.21)
$$\leq 2 \sup_{\widetilde{x} \in D_{i,j}} E_{\widetilde{x}}[|\pi_{\mathbb{Z}}(X_{[0, T_{(B(g_i L, L \log L) \cap A) \times \mathbb{Z}}]} \cap D_{i,j})|]$$

$$\overset{(5.3)}{\leq} cL^{\beta/2}(\log L)^{-\mu - \beta/2}.$$

Proceeding as in (5.16) and (5.19), we thus obtain

$$(5.22) \quad \sup_{z \in E} E_z \left[ \sum_{(i,j) \in \mathcal{H}_0} |\pi_{\mathbb{Z}}(X_{[0, T_{\widetilde{B}'}]} \cap D_{i,j})| \right] \leq cL^{\beta/2}(\log L)^{-\mu - \beta/2} M^{1-\nu},$$

and a similar inequality holds with $\mathcal{H}_1$ in place of $\mathcal{H}_0$.

Using once again a variation on Khasminskii's lemma (cf. (2.46) of [8]), (5.20) and (5.22) imply that for $m = 0, 1$,

$$(5.23) \quad \sup_{z \in E} E_z \left[ \exp \left\{ c \frac{(\log L)^{\mu}}{L M^{1-\nu}} \sum_{(i,j) \in \mathcal{H}_m} |\pi_G(X_{[0, T_{\widetilde{B}'}]} \cap D_{i,j})| \right\} \right] \leq 2$$

and

$$(5.24) \quad \sup_{z \in E} E_z \left[ \exp \left\{ c \frac{(\log L)^{\mu + \beta/2}}{L^{\beta/2} M^{1-\nu}} \sum_{(i,j) \in \mathcal{H}_m} |\pi_{\mathbb{Z}}(X_{[0, T_{\widetilde{B}'}]} \cap D_{i,j})| \right\} \right] \leq 2.$$

We now return to our main objective, that is, bounding the probability in (5.11). We thus see with (5.13), (5.14) and the above controls that, for large $G \in \widetilde{\mathcal{G}}$ and $x \in E$,

$$P_x \left[ S \text{ is not good in } D \text{ and } \sum_{k \geq 1} 1\{\widetilde{R}_k \leq |G|^{2(1-\delta)}\} \leq K_0 \right]$$

(5.25)
$$\leq 2 \exp\{c(\log |G|)^{\rho} \log \log |G|\}$$

$$\times \sup_{\mathcal{H}_0, \mathcal{H}_1} \sup_{m \in \{0,1\}} P_x \left[ \sum_{(i,j) \in \mathcal{H}_m} \sum_{1 \leq k \leq K_0} |\pi_G(X_{[0, T_{\widetilde{B}'}]} \cap D_{i,j})| \right.$$

$$\left. \circ \theta_{\widetilde{R}_k} \geq \frac{L}{2} \left[ \frac{M}{2} \right] \right]$$

$$+ P_x \left[ \sum_{(i,j) \in \mathcal{H}_m} \sum_{1 \leq k \leq K_0} |\pi_{\mathbb{Z}}(X_{[0, T_{\widetilde{B}'}]} \cap D_{i,j})| \circ \theta_{\widetilde{R}_k} \right.$$

$$\geq \frac{1}{2}[L^{\beta/2}(\log L)^{-\beta/2}]\left[\frac{M}{2}\right]\Big]$$

$$\stackrel{(5.23),(5.24)}{\leq} 2\exp\{c(\log|G|)^\rho \log\log|G|\}$$
$$\times 2^{K_0}\left(\exp\left\{-c\frac{LM}{LM^{1-\nu}}(\log L)^\mu\right\}\right.$$
$$\left.+\exp\left\{-c\frac{L^{\beta/2}(\log L)^{-\beta/2}M}{L^{\beta/2}M^{1-\nu}}(\log L)^{\mu+\beta/2}\right\}\right)$$
$$\stackrel{(5.7),(5.9)}{\leq} c\exp\{-c(\log L)^\mu M^\nu\} \stackrel{(5.7)-(5.9)}{\leq} o(|G|^{-2})$$

[with the last inequality of (5.7) and (5.9), $\log L$ is comparable to $\log|G|$]. This shows (5.11) and thus concludes the proof of Theorem 5.2.  □

We now provide an application of Theorem 5.2 in the spirit of Corollary 4.6. We consider an infinite connected graph $G_\infty$ with degree bounded by $d_0 \geq 2$, satisfying the heat kernel bounds of (0.5) [and hence, (0.6) for suitable $\widetilde{\kappa}_i$, $i = 1, 2$], but unlike (4.24), we now assume that

$$\beta \geq 2, \beta \geq \alpha \geq \left(1+\frac{\beta}{2}\right)\vee(\beta-1) \qquad \text{(and therefore, } G_\infty \text{ is recurrent).}$$
(5.26)

We refer to Barlow [3] for examples of such $G_\infty$, when $\beta > 2$, the case $\beta = 2$ being more common. We assume that we have a sequence of finite connected graphs $G_N$, $N \geq 1$, with degree bounded by $d_0 \geq 2$, and $\lim|G_N| = \infty$, and (4.25) and (4.26) hold. We then have following:

COROLLARY 5.3. *Under the above assumptions for all $\delta > 0$, $\varepsilon > 0$,*

(5.27) $$\lim_{N\to\infty}\inf_{x\in E_N} P_x[|G_N|^{2(1-\delta)} \leq \mathcal{T}_{G_N} \leq |G_N|^2(\log|G_N|)^{4+\varepsilon}] = 1.$$

PROOF. We only need to discuss the lower bound on $\mathcal{T}_{G_N}$. We choose $0 < \delta < \frac{1}{2} \wedge \eta$; see (4.26). For the remainder of the proof, all constants $c$ may depend on $d_0, \delta, \alpha, \beta, \kappa_i, 1 \leq i \leq 4$ [cf. (0.5)] and $\eta$. We choose for large $N$ [cf. (4.25) and (5.1)]

(5.28) $$A_N = B(g_N, r_N) \quad \text{and} \quad V_N \text{ a geodesic segment initiating at } g_N \text{ with}$$
$$|V_N| = \left[\frac{r_N}{8} \wedge |G_N|^{\delta/8}\right].$$

With similar arguments as for the proof of (4.29), with the only difference that in the last expression of (4.31), and in the first term in the last member of (4.34), there is an additional factor $T$, since (5.26) replaces (4.24), one has

(5.29) $\quad$ for large $N$, $\quad G_N \in \mathcal{G}_0\Big(d_0, \delta, \delta' = \dfrac{\delta}{8} \wedge \dfrac{\eta}{2\alpha}, \gamma = \dfrac{\alpha}{\beta}, a = \kappa_1\Big).$

In view of Theorem 5.2, our claim will thus follow if we can see that we can select the remaining parameters $\mu, \nu, a', a_G, a_{\mathbb{Z}}$, so that (5.4) holds and

(5.30) $\quad$ for large $N$, $\quad G_N \in \widetilde{\mathcal{G}}_0.$

To this end, we consider the walk on $G_\infty \times \mathbb{Z} \stackrel{\text{def}}{=} E_\infty$, and its Green function

(5.31) $$g_\infty(x, x') = \sum_{k \geq 0} P_x^{E_\infty}[X_k = x'], \qquad x, x' \in E_\infty.$$

LEMMA 5.4 ($\alpha + 1 \geq \beta \geq 2$, $\alpha > 1$). *For a suitable $c > 1$, for any $x, x'$ in $E_\infty$,*

(5.32) $$\frac{1}{c} D(x, x')^{-(2\alpha/\beta - 1)} \leq g_\infty(x, x') \leq c D(x, x')^{-(2\alpha/\beta - 1)} \qquad \text{with}$$
$$D(x, x') = \max(d_{G_\infty}(g, g')^{\beta/2}, |u' - u|, 1),$$
$$x = (g, u), \qquad x' = (g', u').$$

We refer to the Appendix for the proof of this lemma. Then consider $W$ a geodesic segment in $G_\infty$ of length $m \geq 2$, and $J$ an interval of $\mathbb{Z}$ with length $[m^{\beta/2}(\log m)^{-\beta/2}]$. Picking $x_0 = (g_0, u_0) \in W \times J$, with (5.32), one has

(5.33) $$\inf\{g_\infty(x, x_0); x \in W \times J\} \geq c m^{-(\alpha - \beta/2)},$$

and it now follows from the fact that $g_\infty(X_{n \wedge H_{W \times J}}, x_0), n \geq 0$, is a martingale under any $P_x^{E_\infty}$, $x = (g, u) \in E_\infty$, that

(5.34) $$P_x^{E_\infty}[H_{W \times J} < \infty] \leq \frac{g_\infty(x, x_0)}{c m^{-(\alpha - \beta/2)}}$$
$$\stackrel{(5.32)}{\leq} c'\bigg[\max\bigg(\frac{d(g, g_0)}{m}, \frac{|u - u_0|^{2/\beta}}{m}\bigg)\bigg]^{-(\alpha - \beta/2)}.$$

This readily implies that (5.2) holds for large $N$, with an adequate choice of $a' = 1 \vee c'$, with $c'$ as in (5.34) and

(5.35) $$\nu = \alpha - \frac{\beta}{2}.$$

We then continue and check (5.3). Let $W, J$ be as above, and define for $g_0 \in W$, $U_{g_0} = \{g_0\} \times J$ and

$$(5.36) \quad \mathcal{U}_{g_0}(z) = \sum_{x \in U_{g_0}} g_\infty(z, x) \quad \text{for } z \in E_\infty, \text{ so that}$$

$$(5.37) \quad \inf_{z \in U_{g_0}} \mathcal{U}_g(z) \stackrel{(5.32)}{\geq} c \sum_{\ell=1}^{|J|} \ell^{-(2\alpha/\beta - 1)} \geq \begin{cases} c \log m, & \text{if } \alpha = \beta, \\ c|J|^{2(1-\alpha/\beta)}, & \text{if } \alpha < \beta. \end{cases}$$

On the other hand, for $z = (g, u) \in W \times J$ we also have with (5.32)

$$(5.38) \quad \begin{aligned} \mathcal{U}_{g_0}(z) &\leq c \sum_{\ell=1}^{|J|} \frac{1}{d(g, g_0)^{\alpha - \beta/2} \vee \ell^{2\alpha/\beta - 1}} \\ &\leq \begin{cases} c\left(1 + \log\left(\frac{m}{d(g, g_0) \vee 1}\right)\right), & \text{if } \alpha = \beta, \\ c\left(\frac{|J|}{(d(g, g_0) \vee 1)^{\alpha - \beta/2}}\right) \wedge |J|^{2(1-\alpha/\beta)}, & \text{if } \alpha < \beta. \end{cases} \end{aligned}$$

With a similar argument as in (5.34), we find that for $z = (g, u) \in W \times J$,

$$(5.39) \quad P_z^{E_\infty}[H_{U_{g_0}} < \infty] \leq \begin{cases} \dfrac{c}{\log m}\left(1 + \log\left(\dfrac{m}{d(g, g_0) \vee 1}\right)\right), & \text{if } \alpha = \beta, \\ c\left(\dfrac{|J|^{2/\beta}}{d(g, g_0) \vee 1}\right)^{\alpha - \beta/2} \wedge 1, & \text{if } \alpha < \beta. \end{cases}$$

Note that, for $x \in E_\infty$,

$$\begin{aligned} E_x^{E_\infty}&[|\pi_G(X_{[0,\infty)} \cap (W \times J))|] \\ &= E_x^{E_\infty}\left[\sum_{g_0 \in W} 1\{H_{U_{g_0}} < \infty\}\right] \\ (5.40) \quad &= E_x^{E_\infty}\left[H_{W \times J} < \infty, E_{X_{H_{W \times J}}}^{E_\infty}\left[\sum_{g_0 \in W} 1\{H_{U_{g_0}} < \infty\}\right]\right] \\ &\stackrel{(5.39)}{\leq} \max_{z = (g,v) \in W \times J} \sum_{g_0 \in W} \frac{c}{\log m}\left(1 + \log\left(\frac{m}{d(g, g_0) \vee 1}\right)\right) \\ &\leq \begin{cases} c\dfrac{m}{\log m}, & \text{if } \alpha = \beta, \\ c\dfrac{m}{\log m} \log \log m, & \text{if } \alpha < \beta, \end{cases} \end{aligned}$$

using the fact that $W$ is a geodesic segment of length $m$, and the last line of (5.39) together with the inequality $\alpha \geq 1 + \frac{\beta}{2}$ and $|J|^{2/\beta} \approx m/\log m$, when $\alpha < \beta$ [note that the $\log \log m$ factor in the last line of (5.40) is only needed when $\alpha = 1 + \frac{\beta}{2}$]. This takes care of the first estimate in (5.3). For the second estimate, we consider $u_0 \in J$, $W_{u_0} = W \times \{u_0\}$ and

$$\mathcal{W}_{u_0}(z) = \sum_{x \in W_{u_0}} g_\infty(z, x) \qquad \text{for } z \in E_\infty. \tag{5.41}$$

Again with (5.32), we see that

$$\inf_{z \in W_{u_0}} \mathcal{W}_{u_0}(z) \geq c \sum_{\ell=1}^{m} \frac{1}{\ell^{\alpha - \beta/2}} \geq \begin{cases} c, & \text{if } \alpha > 1 + \frac{\beta}{2}, \\ c \log m, & \text{if } \alpha = 1 + \frac{\beta}{2}. \end{cases} \tag{5.42}$$

When $z = (g, u) \in W \times J$, we also find

$$\mathcal{W}_{u_0}(z) \leq c \sum_{\ell=1}^{m} \frac{1}{\ell^{\alpha - \beta/2} \vee |u - u_0|^{2\alpha/\beta - 1}}$$

$$\leq \begin{cases} c(|u - u_0| \vee 1)^{2/\beta(1 + \beta/2 - \alpha)}, & \text{if } \alpha > 1 + \frac{\beta}{2}, \\ c\left(1 + \log\left(\frac{m}{|u - u_0|^{2/\beta} \vee 1}\right)\right), & \text{if } \alpha = 1 + \frac{\beta}{2}. \end{cases} \tag{5.43}$$

As in (5.39), we find that, for $z = (g, u) \in W \times J$,

$$P_z^{E_\infty}[W_{u_0} < \infty]$$

$$\leq \begin{cases} c(|u - u_0| \vee 1)^{2/\beta(1 + \beta/2 - \alpha)}, & \text{if } \alpha > 1 + \frac{\beta}{2}, \\ \frac{c}{\log m}\left(1 + \log\left(\frac{m}{|u - u_0|^{2/\beta} \vee 1}\right)\right), & \text{if } \alpha = 1 + \frac{\beta}{2}. \end{cases} \tag{5.44}$$

A similar computation as in (5.40) shows that, for $x \in E_\infty$,

$$P_x^{E_\infty}[|\pi_{\mathbb{Z}}(X_{[0,\infty)} \cap (W \times J))|]$$

$$\leq \max_{z=(g,u) \in W \times J} \sum_{u_0 \in J} c(|u - u_0| \vee 1)^{2/\beta(1+\beta/2-\alpha)}$$

$$\leq \begin{cases} cm^{\beta/2 + (1+\beta/2-\alpha)}, & \text{if } \alpha > 1 + \frac{\beta}{2}, \\ c\frac{m^{\beta/2}}{(\log m)^{\beta/2+1}} \log \log m, & \text{if } \alpha = 1 + \frac{\beta}{2}. \end{cases} \tag{5.45}$$

Combining (5.40) and (5.45), we see that with a suitable choice of $a_G, a_{\mathbb{Z}} \geq 1$, and

$$\mu < 1, \tag{5.46}$$

condition (5.3) holds for large $N$. This concludes the proof of (5.30), with parameters [cf. (5.35) and (5.46)] that fulfill (5.4). □

In the same spirit as Remark 4.7, we have the following:

REMARK 5.5. (i) Corollary 5.3 applies, in particular, when $G_N$ is some ball of radius $N$ in $G_\infty$ (with possibly variable center), with $G_\infty$ as above in Corollary 5.3.

(ii) We can choose $G_\infty$ to be $\mathbb{Z}^2$, with its usual graph structure, corresponding to $\alpha = \beta = 2$ in (5.26). Corollary 5.3 then applies to the case $G_N = (\mathbb{Z}/N\mathbb{Z})^2$, $N \geq 1$. This with Remark 4.7(ii) recovers Theorems 1.1 and 2.1 of [8], for arbitrary $d \geq 1$. Of course, Corollary 5.3 applies just as well to $G_N = [0, N]^2$, $N \geq 1$.

(iii) When $G_N = [0, [N^\lambda]] \times [0, N]$, $N \geq 1$, with $\lambda \in (0, \infty)$, then for large $N$, $\lambda_{G_N} N^{2(\lambda \vee 1)}$ remains bounded and bounded away from zero. Picking $A_N$ as a ball with suitable center and radius a small multiple of $N^{\lambda \wedge 1}$, we can apply Corollary 5.3, as soon as [cf. (4.26) and Remark 4.7(iii)], $\frac{1}{2} < \lambda < 2$.

## APPENDIX

We prove here Lemma 5.4. We can clearly replace $g(\cdot, \cdot)$ in (5.32) with

$$(A.1) \qquad \overline{g}_\infty(x, x') = E_x^{E_\infty}\left[\int_0^\infty 1\{\overline{X}_t = x'\}\,dt\right], \qquad x, x' \in E_\infty,$$

since for a suitable constant $c \geq 1$, $\frac{1}{c} \leq \overline{g}_\infty/g_\infty \leq c$. We then use a similar representation as in Lemma 3.3. We introduce the continuous (piecewise linear) increasing processes

$$(A.2) \qquad \begin{aligned} A_t &= \int_0^t \deg(\overline{Y}_r)\,dr, \qquad t \geq 0, \quad \text{and} \\ \tau_s &= (A^{-1})_s = \inf\{t > 0; A_t > s\}, \qquad s \geq 0, \end{aligned}$$

and the time changed process

$$(A.3) \qquad \widetilde{Y}_s = \overline{Y}_{\tau_s},$$

so that under $P_g^{G_\infty}$, $\widetilde{Y}.$ is the continuous time walk on $G_\infty$ with constant jump rate equal to 1, starting at $g \in G_\infty$. We also note that with (A.2) and (A.3)

$$(A.4) \qquad \tau_s = \int_0^s \deg(\widetilde{Y}_r)^{-1}\,dr, \qquad s \geq 0.$$

With similar calculations as in Lemma 3.3, we thus find that, for $x = (g, u)$, $x' = (g', u')$ in $E_\infty$, one has

$$(A.5) \qquad \overline{g}_\infty(x, x') = E_g^{G_\infty} \otimes E_u^\mathbb{Z}\left[\int_0^\infty 1\{\widetilde{Y}_s = g'\}1\{\overline{Z}_{\tau_s} = u'\}\deg(\widetilde{Y}_s)^{-1}\,ds\right],$$

and therefore, for a suitable constant $c \geq 1$,

$$\frac{1}{c} \leq \frac{1}{g_\infty(x,x')} \int_0^\infty E_g^{G_\infty}[\widetilde{Y}_s = g', q(\tau_s, u, u')]\, ds \leq c$$

(A.6)

with $q(t, u, u') = P_u^{\mathbb{Z}}[\overline{Z}_t = u']$.

We first show the right-hand side inequality of (5.32). As in (3.21), we denote with $N_t, t \geq 0$, a Poisson counting process with unit intensity. We set $t_0 = D$, with $D$ as in (5.32) and note that with (3.21), $P[N_t \geq 2D] \leq P[N_{t_0} \geq 2D] \leq ce^{-cD}$, for $t \leq t_0$, so that

$$\int_0^\infty E_g^{G_\infty}[\widetilde{Y}_s = g', q(\tau_s, u, u')]\, ds$$

$$\leq t_0 \left( ce^{-cD} + \sup_{k \leq 2D} P_g^{G_\infty}[Y_k = g'] \sup_{\ell \leq cD} P_u^{\mathbb{Z}}[Z_\ell = u'] \right)$$

$$+ \int_{t_0}^\infty \left( ce^{-ct} + \sup_{t/2 \leq k \leq 2t} P_g^{G_\infty}[Y_k = g'] \sup_{ct \leq \ell \leq c't} P_u^{\mathbb{Z}}[Z_\ell = u'] \right) dt$$

$$\overset{(0.5)(i)}{\leq} t_0 \left( ce^{-cD} + c\exp\left\{ -c\left(\frac{d_{G_\infty}(g,g')^\beta}{D}\right)^{1/(\beta-1)} - c\frac{|u-u'|^2}{D} \right\} \right)$$

(A.7)

$$+ ce^{-ct_0} + \int_{t_0}^\infty \frac{c}{t^{\alpha/\beta+1/2}}$$

$$\times \exp\left\{ -c\left(\frac{d_{G_\infty}(g,g')^\beta}{t}\right)^{1(\beta-1)} -c\frac{|u-u'|^2}{t} \right\} dt$$

$$\leq cDe^{-cD^{1/(\beta-1)}} + \frac{c}{D^{2\alpha/\beta-1}} \leq \frac{c}{D^{2\alpha/\beta-1}},$$

using the change of variable $t = D^2 s$ to bound the last integral. With (A.6), this proves the right-hand side inequality of (5.32). To prove the left-hand side inequality of (5.32), we note in addition to (3.21) that, for a suitable constant $c \geq 1$,

$$\frac{1}{c} \leq P[N_t = k+1]/P[N_t = k] \leq c$$

(A.8)

for all $k \in \left[\frac{t}{2}, 2t\right]$ and $t \geq 1$.

We now write with (0.5)(ii) and a similar bound for the walk on $\mathbb{Z}$:

$$g_\infty(x,x') \overset{(A.6)}{\geq} c \int_{D^2}^\infty E_g^{G_\infty}[\widetilde{Y}_s = g', q(\tau_s, u, u')]\, ds$$

(A.9)
$$\geq \int_{cD^2}^{\infty} \frac{c}{s^{\alpha/\beta+1/2}} \, ds = \frac{c}{D^{2\alpha/\beta-1}}.$$

This proves the left-hand side inequality of (5.32) and concludes the proof of Lemma 5.4.

**Acknowledgment.** We wish to thank Amir Dembo for many helpful discussions.

DEPARTEMENT MATHEMATIK
ETH ZURICH
CH-8092 ZÜRICH
SWITZERLAND
E-MAIL: sznitman@math.ethz.ch